\documentclass[onecolumn,showpacs,aps,amsmath,amssymb,superscriptaddress,nofootinbib,11pt]{revtex4-1}

\usepackage{graphicx}
\usepackage{amssymb}
\usepackage{amsmath}
\usepackage{epstopdf}

\usepackage{epigraph}
\usepackage{comment}
\usepackage{color}
\usepackage[breakable]{tcolorbox}
\newcommand{\bea}{\begin{eqnarray}}
\newcommand{\eea}{\end{eqnarray}}

\begin{document}
\title{Forest expansion of two-body partition functions for sparse interaction graphs}
\author{F. Caravelli}
\affiliation{Theoretical Division (T4) and Center for Nonlinear Studies,\\
Los Alamos National Laboratory, Los Alamos, New Mexico 87545, USA}
%\date{Mar 2020}
\begin{abstract}
We study tree approximations to classical two-body partition functions on sparse and loopy graphs via the Brydges-Kennedy-Abdessalam-Rivasseau forest expansion. We show that for sparse graphs (with large cycles), the partition function above a certain temperature $T^*$ can be approximated by a graph polynomial expansion over forests of the interaction graph. Within this ``forest phase", we show that the approximation can be written in terms of a reference tree $\mathcal T$ on the interaction graph, with corrections due to cycles. From this point of view, this implies that high-temperature models are easy to solve on sparse graphs, as one can evaluate the partition function using belief propagation. We also show that there exist a  high- and low-temperature regime, in which $\mathcal T$ can be obtained via a maximal spanning tree algorithm on a (given) weighted graph. We study the algebra of these corrections and provide first- and second-order approximation to the tree Ansatz, and give explicit examples for the first-order approximation. 
\end{abstract}

\maketitle
%\tableofcontents
\section{Introduction}
There has been large interest in the study of statistical models on arbitrary graphs ever since Bethe introduced the notion of ``superlattices" \cite{bethe,baxter}. 
There are multiple reasons for such interest, and for the purpose of this paper we focus on the computability of the partition function for two-body Hamiltonians; in particular, we are interested in tree approximations to the interaction pattern.
The partition function for a statistical model can be interpreted as the generator of moments of the underlying fluctuating variables. The key issue is that the partition function is often hard to calculate. In some cases this problem falls in complexity class from a NP-Hard problem to a polynomial problem. An example of such reduction is for instance in the case of the Ising model on trees. It is known in fact that via cavity method (belief propagation, message passing) \cite{cavity,Montanari,Pearl1,Pearl2,yedidia} it is possible to solve for the partition function or the probability distribution in polynomial time and resources, both for regular and disordered statistical systems. However, proofs of convergence of the Belief Propagation (BP) algorithm are in general strongly restricted to graphs with (locally) tree-like topology, e.g. exact trees or graphs with large cycles \cite{dembo}.
In the case in which the cycle space is non-trivial, work done by Frey and MacKay has shown however that Belief Propagation can sometime converge also in the case in which the graph contains cycles \cite{mackay}. A cycle is a sequence of edges in a graph $\mathcal G$ such that starting from one vertex, one can return to the starting point without passing from the same edge. In the last decade, first Chertkov and Chernyak \cite{chertkov,chertkov2}  have introduced the notion of cycle calculus in belief propagation, showing that for graphs in which cycles are large, BP should be expected to converge.  More recently such cycle expansions was also discussed by Cantwell and Newman in the context of message passing \cite{newman}. 

The result of the expansion is similar to the Mayer or polymeric expansion of a partition function, as it provides a sum over terms of a graph polynomial whose variables are $x_{ij}=e^{\beta A_{ij} H_{ij}}-1$, thus apt for a high temperature expansion. As a matter of fact, the proof is done at the level of the probability distribution, thus without relying on belief propagation. In this sense, our results are similar in spirit to the polymeric expansion approach by Vuffray and Macris \cite{vuffray,vuffray2}, introduced in order to understand ``loopy" belief propagation (when standard belief propagation algorithms converge on loopy graphs) \footnote{In this paper we will use the nomenclature ``cycle" for loop.}. 
If $Y_i$ represents groups of interacting variables, in the statistical mechanics literature,  a polymeric expansion of a partition function takes the form:
\begin{eqnarray}
    \mathcal Z=\sum_{p}\frac{1}{p!} \sum_{\text{disjoint polymers }Y_i} z(Y_i)
\end{eqnarray}
where the $z(Y_i )'s$ are called the activity of
the the polymer $Y_i$, and is also  called a cluster expansion for the polymer gas
while in the (constructive) field theory
literature it is rather called the Mayer expansion. Such expansion can also be used, for instance in the case of the Ising model, to construct ``continuum"  field theories of the Ising model \cite{caravellicont}.

In this paper we obtain, via an alternative derivation which takes advantage of the forest expansion (or BKAR formula)  \cite{brydges,rivaabd}, the loop or cycle expansion previously obtained, clarifying that the loop expansion is contained in a ``forest phase". In fact, the BKAR formula we use in this paper can be considered as a generalization of the polymeric expansion. The paper is organized as follows. In Sec. II we provide a list of the results obtained in this paper, with the minimum amount of technical jargon for clarity, and  provide an interpretation. In Sec. III and IV  we obtain the main results of this paper in a technical fashion. Conclusions follow. 
\section{Main Results}\label{sec:results}
In this paper we obtain, via an independent approach, a cycle expansion for the underlying partition function.  We focus on two-body Hamiltonians for a classical system. Such flexibility is due to the fact that the character of the expansion is focused on the interaction pattern rather than the nature of operators, and this allows us to focus on topology of the system rather than the details of the two body Hamiltonian. The only requirement is that the graph is \textit{sparse}, e.g. that the minimum cycle length grows with the size of the system.  A list of the result that are in this paper, and of which we provide an interpretation below, are contained in Sec. II and III.  In Sec. III.  we introduce the technical results, in particular the application of the BKAR  (IIIA)formula to two-body partition functions. In Sec. IIIB and IIIC we provide the proofs of the forest phase. In Sec. IV we show that the loop expansion is contained, perturbatively, within the forest phase, and we show that this can be obtained via a series in the forms of subtrees.

The partition function, once we fix the nature of the interactions (for instance, the Ising model with classical spins via exchange coupling), depends on the topological aspects of the interactions.

\begin{figure*}
   % \centering
    \includegraphics[scale=1.5]{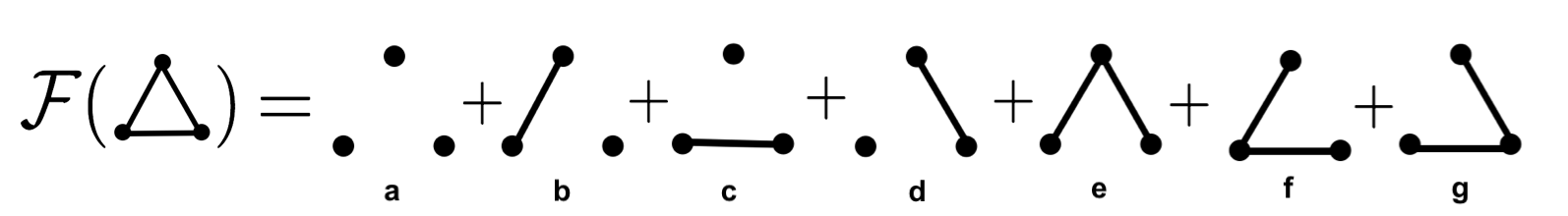}
    \caption{Forest expansion for $K_3$, the complete graph on 3-nodes. }
    \label{fig:forestk3}
\end{figure*}

First of all, let us clarify what are the main ingredients in the partition function approximation, at high temperatures, which applies to sparse graphs. Let $\mathcal G$ be interaction graph. Then we call $\mathcal F(\mathcal G)$ as the sum over all possible connected or disconnected subtrees of $\mathcal G$. An example of such expansion is presented in Fig. \ref{fig:forestk3}.

Specifically, in this paper we provide a proof of various results.
\begin{itemize}
\item The \textit{first result} we obtain in this paper is, for tree-like graph as we discuss below, what has to be considered high and low temperature.  Let $H=\sum_{(ij)\in \mathcal G} H_{ij}=\sum_{ij} A_{ij} H_{ij}$ be a two body Hamiltonian on a graph $\mathcal G$ with adjacency matrix $A_{ij}$, and let the partition function be weighted by the Boltzmann factor $e^{-\beta H}$. We assume that $H$ is a classical Hamiltonian. Let us define $L^*$ to be the minimum length among the cycles in $\mathcal G$. Then, given 
 \begin{eqnarray}
     \beta \leq W(1)\ \text{sup}{|H_{ij}|}=\beta^*
 \end{eqnarray}
 where $W(x)$ is the Lambert-W function (which is the solution of $x e^x=y$), we have
\begin{eqnarray}
    \lim_{L^*\rightarrow \infty} \mathcal Z_{\mathcal G}(\beta<\beta^*)=Z_{\mathcal F}(\beta)
\end{eqnarray}
which is a partition function defined via a forest sum on the interaction graph, e.g. connected or disconnected trees. In the formula above, $Z_{\mathcal F}$ represents the sum of partition functions defined on trees of the forest expansion, and thus without cycles. This implies that every single term can be in principle computed efficiently. From the definition above we see that a requirement is the one of bounded interactions. The result above establishes that for these graph there exist a specific temperature between a ``high" and a ``low" temperature regime.
We will make more precise in what sense a graph is dense or sparse in Section \ref{sec:expansion} where the derivation is presented.
\end{itemize}
Let us clarify the implications and the meaning of the statement above, and also briefly explain the strategy of the proof. The BKAR expansion is a forest expansion as the one in the example of Fig. \ref{fig:forestk3}. By construction of the formula, whose details are presented below, is such that for every present edge there is a variable $u_{ij}$ which is being integrated. If the underlying graph over which the expansion is performed does contain a cycle including a sequence of edges, then the formula includes terms of the form $\text{min }u_{ij}$, where $u_{ij}$ are the variables associated with the existing edges of the particular graph of the expansion. This implies that the forest expansion is not purely a forest (trees) expansion, but implicitly includes loops. The first statement above however guarantees that, in the high temperature phase and if the loops are large, every single term containing a loop goes to zero.  As a result, the expansion implies that if loops are large can be neglected in the high temperature phase. This might seem remarkable at first, since we know that every single graph could be evaluated (depending on the model) efficiently. However, the number of such graph is typically large. Thus, since we are using the limit $\lim_{L\rightarrow \infty} Z_{\mathcal F}$, such result is essentially only a partial achievement, which will be however expanded thanks to the second statement below, which is concerned with the loop algebra. Essentially, we are claiming that for a generic two-body classical partition function, the partition function can be written as $\mathcal Z=Z_{\mathcal F}+\mathcal O(L^*)$, where $L^*$ is the minimum loop length of the graph, with $\lim_{L^*\rightarrow \infty} \mathcal O(L^*)=0$ in the high temperature regime.  A comment that might be interesting to make at this stage is how fast these terms go to zero. What we prove is that $\lim_{L^*\rightarrow } |\mathcal Z-\mathcal Z_{\mathcal F}|=\lim_{L^*\rightarrow \infty} \mathcal O(L^*)$ in the high temperature phase.  We are thus assuming a very specific order of the limits, in which a graph might have large loops, and thus we are taking the limit $N\rightarrow \infty$ together with the limit $L^*\rightarrow \infty$. These can thus be interpreted as tree-like graphs, and this is the sense of statement above. Regarding the difference $\mathcal O$, between the partition function and the forest approximation, we show that every term in $\mathcal O(L^*)$ is a sum of terms which go exponentially to zero, at most as large as
\begin{eqnarray}
    \mathcal O(L^*)=\begin{cases}
    O(q^{L^*-1} e^{\rho \beta \bar H}) & \text{cycle dense graphs} \\
   O( \frac{q^{L^*-1}}{L^*}) & \text{cycle sparse graphs}
    \end{cases}\nonumber
\end{eqnarray}
with $q =\beta W(1) \bar H$ and $\bar H=\text{sup}_{ij}|H_{ij}|$. The definition of cycle dense or cycle sparse and $\rho$ is technical and provided in the bulk of the paper, but it should clarify where the temperature $\beta^*$ originates from. Note that the fact that the forest expansion is related to a forest, it does not in fact mean that the partition function is the one of a tree. This is clarified in the result below.

\begin{itemize}
\item The \textit{second result} is given by
\begin{eqnarray}
    \lim_{\beta \rightarrow 0^+}\mathcal Z_{\mathcal F}(\beta) =\mathcal Z_{\mathcal T}+O(\beta \bar H)
\end{eqnarray}
where $\mathcal Z_{\mathcal T}$ is defined on an arbitrary spanning tree of the original interaction $\mathcal G$. The rest $\mathcal R$ depends on the properties of the cycles, which we will discuss below.
\end{itemize}

\begin{figure}
    \centering
    \includegraphics[scale=1.5]{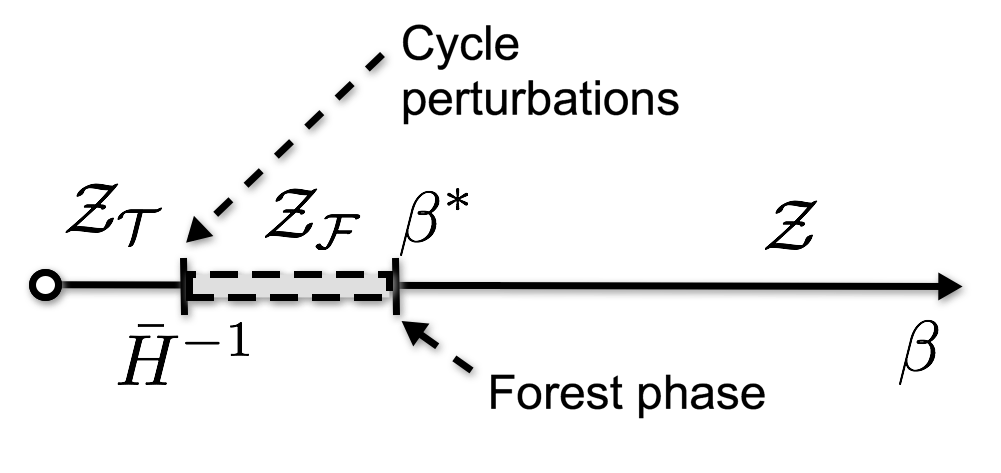}
    \caption{Temperature hierarchy derived ``sparse" graphs: the ``cycle perturbative" regime is hidden in a ``forest" phase. }
    \label{fig:cyclepert}
\end{figure}

The notation above can be clarified via a simple example. Consider for instance Fig. \ref{fig:forestk3}. The first result shows that (if the loop would be large and not only of size 3), that we could approximate the partition function as a sum over the tree subgraphs. However, it is useful to actually identify a subset of graphs which would be associated to a tree. Consider for instance, again in Fig. \ref{fig:forestk3}, the graphs \textit{a,b,c,f}. These is the forest of a tree which contain only two of the edges of the original interaction graph, and it can be shown that we can re-sum them to be the partition function of the system on a certain reference tree. Clearly, this leaves graphs $d,e,g$ out. As we show in the paper, the difference between $\mathcal Z_{\mathcal F(\mathcal G)}$ and $\mathcal Z_{\mathcal T}$ is of order $O(\beta \bar H)$ (again, assuming $\lim_{L^*\rightarrow \infty} \mathcal Z$). The situation is thus the one presented in Fig. \ref{fig:cyclepert}. There is an intermediate temperature region between $\beta^*$ and $\beta_{HT}=\frac{1}{\bar H}$ in which we can write $\mathcal Z_{\mathcal F}$ in terms of a reference tree. Such perturbative regime is described by the cycle algebra obtained in Sec. \ref{sec:cyclealg}, and is also the situation analyzed by \cite{chertkov} and \cite{vuffray}. Since in principle, referring again to Fig. \ref{fig:forestk3}, we could have chosen 2 other subtrees, the question is which of these subtrees could have been chosen. The answer is simple, and it is the tree for which the couplings not in the reference tree are the smallest.

Thus, the strategy of the analysis is in two steps. In the first step we show that for temperature such that $\beta\leq \beta^*$ we can approximate the partition function using $\mathcal Z_{\mathcal F}$, the forest expansion. In the second step we show that $\mathcal Z_{\mathcal F}\approx \mathcal Z_{\mathcal T}+O(\beta \bar H)$. The result thus parallels the one by Chertkov and Chernyak, with the caveat that we highlight this intermediate phase.

The first statement is derived using the Brydges-Kennedy-Abdessalam-Rivasseau (BKAR) forest formula \cite{brydges,rivaabd,gurau} for functions of $N(N-1)/2$ variables, and thus for a finite size system. 
In our case these $N(N-1)/2$ variables are associated with the interaction graph $A_{ij}$. The advantage of such an approach is that it provides an alternative to belief propagation in deriving cycle expansions. While the expansion is based on a forest (e.g. sum over connected or disconnected trees), cycles enters subtly into BKAR forest expansion. We show that if the length $L$ of a cycle is large and $\beta<\beta^*$, terms of the exact expansion for the probability distribution containing the cycles are suppressed by a factor $\frac{q^{L-1}}{L}$ where $q<1$. As a result, this implies that above a certain temperature we can approximate our probability distribution with a tree (plus cycle corrections).

Regarding the first statement, interestingly enough this approach gives an evaluation of the function in terms of forests on the interaction graph $\mathcal G$ and its adjacency matrix $A_{ij}$; the sum it is over disconnected tree expansion, where however cycles enter in a non-trivial manner.

The second statement is instead graph theoretical in nature. We show that the forest expansion for a generic graph written in terms of a reference spanning tree. In particular, we show that the introduction of higher order cycles has a hierarchy in temperature, e.g. they contribute more and more as the temperature is lowered.

It is also important to notice that the forest phase can in principle provide a way to evaluate the partition function numerically. In fact, every single term of the forest expansion is a tree, which can be evaluated in polynomial time. Thus, performing a Monte Carlo on random (connected) subtrees it might be possible to evaluate $\log Z$.

%We then show that the forest sum can be written in terms of a reference tree, and establish that there exist an optimal tree which can be obtained from a maximal spanning tree on a graph whose weights depend on the couplings and the temperature.
%We then provide corrections to the tree approximation in the high temperature expansion. %Conclusions follows.

\section{Forest expansion for two-body partition functions}\label{sec:expansion}

We consider the general partition function
\begin{eqnarray}
     \mathcal Z=\int D\mu(\sigma) e^{\beta \sum_{ij} A_{ij} H_{ij}(\sigma_i,\sigma_j)},
\end{eqnarray}
e.g. $H=-\sum_{ij} A_{ij} H_{ij}$, for an arbitrary measure $D\mu(\sigma)$, but for reasons that will become clear soon we restrict our approach to two body interactions. Given this setup,  here we first propose to rewrite the partition function as a forest expansion, and in particular a reduction in terms of trees. 
Let us define $x_{ij}=e^{\beta A_{ij} H_{ij}}-1$.
A forest expansion is a sum of the form
\begin{eqnarray}
    \mathcal Z=\sum_{\mathcal F\in\mathcal G } Z_{\mathcal F}(x_{ij})
\end{eqnarray}
in which $\mathcal F$ are subgraphs of $\mathcal G$ which do not contain cycles. Clearly, the fact that $\mathcal F$ does contain cycles does not mean that cycles effectively do not contribute to the sum, and we will see how these enter in the final expression in a perturbative approach (high temperature).

 In the following we will often refer to $\sup |\beta H_{ij}|'s$ or  $\sup |e^{\sum_{ij} u_{ij}\beta H_{ij}}|'s$.%\textcolor{red}{check that this is true in all the proofs}

Since the partition function depends on the interaction graph topology,
, we write  the following explicit expression for the partition function
\begin{equation}
    \mathcal Z(\beta,\{A_{ij}\})=Z(\beta,A_{12},A_{23},\cdots)
\end{equation}
where we consider the variables $A_{ij}$ representing the couplings as $N(N-1)/2$ variables, and introduce dummy variables $h_{ij}$ as
\begin{equation}
   f(h_{ij})=\mathcal Z(\beta,A_{12}h_{12},A_{23}h_{23},\cdots)
\end{equation}
It follows that below we will interpret $\mathcal Z$ to be explicitly a function of the interaction graph $\mathcal G$.
Given this approach, we now introduce the BKAR formula for this representation. The focus is to calculate $\mathcal Z(\beta)=f(\vec 1)$.

\subsection{Brydges-Kennedy-Abdessalam-Rivasseau Forest formula}
The next step of our approach is the use of an exact expression for the partition function in terms of trees. We thus employ the Brydges-Kennedy-Abdessalam-Rivasseau (BKAR) formula in order to write $f(\vec 1)$ in terms of a forest on the interaction graph \cite{brydges,rivaabd,gurau}.

In order to make the notation clear, we consider $N$ points labeled by $i=1,\cdots,N$ which represent the number of variables $\sigma$ over which the partition function is summed or integrated over, and denote $(ij)$ with edges between nodes $i$ and $j$. We then consider a smooth and arbitrarily derivable function $f([0,1]^{N(N-1)/2})\rightarrow \mathbb{R}$. If the system (and thus $N$) is finite this is the case, but in the limit $N\rightarrow \infty$ such ability to take derivatives can be undermined by phase transitions.  While restrictive, there are many other techniques to evaluate partition functions close to phase transitions, and thus we focus on the case in which $N$ is large but finite. We wish to employ the following forest formula:
\begin{equation}
    f(1,\cdots,1)=\sum_{\mathcal F} \underbrace{\int_0^1\cdots \int_0^1}_{|\mathcal F|}  \Big( \big(\prod_{(ij)\in \mathcal F} \frac{\partial}{\partial u_{ij} }\big) f\Big)|_{u_{ij}=w_{ij}(u_{\mathcal F})}.
\end{equation}
where the sum runs over all forests $\mathcal F$ drawn over the labeled vertices. While the notation of the formula above is a bit intimidating, one essentially takes the derivative of the function with respect of the variables of the forest elements, and then sets the variables to $w_{ij}$ after the derivative is taken; the function, which we define below,  depends on the structure of the forest and the variables $u_{ij}$.

A forest is a sum over trees (not necessarily spanning the graph), connected or disconnected, over the complete graph $K_{N}$ on the $N$ vertices. The typical example is the one in Fig. \ref{fig:forestk3}. The derivatives in the forest expansion are evaluated at the point $u_{ij}=w_{ij}(\mathcal F)$. Given a connected tree $\mathcal T$ which contains vertices $i$,$j$ we call $\mathcal P_{ij}$ the unique sequence of edges on the tree $\mathcal T$, $\mathcal P_{ij}^{\mathcal T}=\{e_{ik},e_{kt},\cdots,e_{rj}\}$. It is important to note that $\mathcal P_{ij}^{\mathcal T}$ is always unique, as every element of the sum is a tree.
Then

\begin{equation}
    w_{ij}(\mathcal F)=\begin{cases}
    u_{ij} & \text{if }(ij)\in \mathcal F\\
    \text{min}(u_{kt}\in P_{ij}^{\mathcal F})& \text{if  } \exists \text{ a nonempty }P_{ij}\\
    0 & \text{otherwise}.
    \end{cases}
\end{equation}
Albeit the formula applies to the complete graph $K_N$ of the interaction matrix,  in practice it is evaluated on the support of $\mathcal G$. This is due to the fact that the BKAR expansion contains terms of the form $\prod_{(ij)\in \mathcal F} \frac{\partial}{\partial u_{ij} } f(\{ u_{ij} A_{ij}\})$ will bring down from the exponential $A_{ij}$, which are zero if the forest includes a coupling term which is absent.

Let us now evaluate the derivatives explicitly. We have
\begin{eqnarray}
    \Big(\prod_{(ij)\in \mathcal F} \frac{\partial}{\partial u_{ij} } \Big)f=\prod_{(ij)\in \mathcal F} \beta A_{ij} H_{ij} f(\{h\})
\end{eqnarray}
which follows from the fact that $f(\{u\})=\int d D\mu(\sigma) e^{\sum_{ij} A_{ij} H_{ij} u_{ij}}$.
Thus, we have
\begin{widetext}
\begin{equation}
    f(1,\cdots,1)=\sum_{\mathcal F} \int d\mu(\sigma) \underbrace{\int_0^1\cdots \int_0^1}_{|\mathcal F|}  \prod_{(ij)\in \mathcal F}   du_{ij} \beta A_{ij} H_{ij} e^{\beta \sum_{ij} A_{ij} H_{ij} w_{ij}}|_{w_{ij}(u_{\mathcal F})}.
\end{equation}
\end{widetext}
Let us now look at the following different cases. 

Let $\mathcal F$ be a forest and $\mathcal G$ the interaction graph. Let us call $\bar {\mathcal G}=\mathcal G\setminus \mathcal F$ the remainder of the graph. We define a path augmentation as the introduction of $(ij)\in \bar{\mathcal G}$ in $\mathcal F$.
If there is no path augmentation, e.g. $(ij) \in \bar {\mathcal G}$ such that $\mathcal F\cup (ij)$ contains a cycle, then $w_{ij}$ has only support on $\mathcal F$, and $w_{ij}=u_{ij}$. Let us call $\mathcal A$ the subset of forests such that there is no path augmentation and $\bar {\mathcal A}=\mathcal F\setminus \mathcal A$.  Clearly, $$\sum_{\mathcal F}=\sum_{\bar {\mathcal  A}}+\sum_{\mathcal A}.$$

In the case in which there is no path augmentation it is not hard to calculate the integrals explicitly, as these do not depend on $\text{min}(u's)$. We have 
\begin{eqnarray}
    \mathcal Z_{\bar A}&=&\int d\mu(\sigma) \int_0^1\cdots \int_0^1  \prod_{(ij)\in \mathcal  A}   du_{ij} \beta A_{ij} H_{ij} e^{\beta \sum_{ij} A_{ij} H_{ij} u_{ij}}  \nonumber \\
    &=&\int d\mu(\sigma)  \prod_{(ij)\in \mathcal  A} (e^{\beta A_{ij} H_{ij}}-1).
\end{eqnarray}
It is thus convenient to rewrite the sum by dividing into terms with and without path augmentations.

Let us write
\begin{widetext}
\begin{eqnarray}
    \mathcal Z &=&\sum_{\bar {\mathcal A}} \int d\mu(\sigma) \underbrace{\int_0^1\cdots \int_0^1}_{|\mathcal F|}  \prod_{(ij)\in \mathcal F}   du_{ij} \beta A_{ij} H_{ij} e^{\beta \sum_{ij} A_{ij} H_{ij} w_{ij}}\nonumber \\
    &+&\sum_{ {\mathcal A}} \int d\mu(\sigma) \underbrace{\int_0^1\cdots \int_0^1}_{|\mathcal F|}  \prod_{(ij)\in \mathcal F}   du_{ij} \beta A_{ij} H_{ij} e^{\beta \sum_{ij} A_{ij} H_{ij} u_{ij}} e^{\beta \sum_{k} A_{i^\prime_k j^\prime_k} H_{i^\prime_k j^\prime_k} \text{min}(u_{i_k j_k} \in \mathcal P_{i^\prime_k j^\prime_k})} \nonumber \\
    &=&\sum_{{\mathcal F}} \int d\mu(\sigma) \underbrace{\int_0^1\cdots \int_0^1}_{|\mathcal F|}  \prod_{(ij)\in \mathcal F}   du_{ij} \beta A_{ij} H_{ij} e^{\beta \sum_{ij} A_{ij} H_{ij} w_{ij}}\nonumber \\
    &+&\sum_{ {\mathcal A}} \int d\mu(\sigma) \underbrace{\int_0^1\cdots \int_0^1}_{|\mathcal F|}  \prod_{(ij)\in \mathcal F}   du_{ij} \beta A_{ij} H_{ij} e^{\beta \sum_{ij} A_{ij} H_{ij} u_{ij}} (e^{\beta \sum_{k} A_{i^\prime_k j^\prime_k} H_{i^\prime_k j^\prime_k} \text{min}(u_{i_k j_k} \in \mathcal P_{i^\prime_k j^\prime_k})}-1) \nonumber \\
    &=&\mathcal Z_{\mathcal F(\mathcal G)}+\mathcal Z_{\mathcal A},
    \label{eq:forestred}
\end{eqnarray}
\end{widetext}
where $\mathcal Z_{\mathcal F(\mathcal G)}$ is a sum in which cycles have been removed; on the other hand $\mathcal Z_{\mathcal A}$ contains all the cycles, which we denoted $\mathcal P_{ij}$, where $(i,j)$ is the edge being added to form a cycle in the forest element. In the following, we will refer to $\mathcal Z_{\mathcal F}$, implicitly referring to the forest expansion derived from the interaction graph $\mathcal G$. In the equation above we have simply added and subtracted one to all the exponentials containing cycle augmentations
 \begin{eqnarray}
     e^{\beta \sum_{k} A_{i^\prime_k j^\prime_k} H_{i^\prime_k j^\prime_k} \text{min}(u_{i_k j_k} \in \mathcal P_{i^\prime_k j^\prime_k})}+1-1.
 \end{eqnarray}
Clearly, if $\mathcal G$ is a tree, then there are no path augmentations, and the formula is a sum over all possible forests. In this case, if $\mathcal G$ is a tree, we have
\begin{eqnarray}
    \mathcal Z=f(\vec 1)=\sum_{\mathcal F} \int d\mu(\sigma) \prod_{(ij) \in \mathcal F}( e^{\beta A_{ij} H_{ij}}-1)=Z_{\mathcal F},\nonumber
\end{eqnarray}
which is exactly the cluster expansion in statistical physics at the tree level.
We will prove below that if the underlying graph $\mathcal G$ is a tree, such sum is exactly $Z_{\mathcal G}=\int d\mu\ e^{-\beta H}$. Thus, this proves that we recover the partition function of the 2-body Hamiltonian on a tree. 

What is interesting is the case in which $\mathcal G$ is not a tree, which by definition implies that it contains cycles. 
Let thus us now consider the most interesting case, e.g. if $\mathcal F\notin \mathcal A$.
Since we want to eventually obtain an expansion in terms of the minimum length $L$, let us assume that the underlying graph $\mathcal G$ has a minimum cycle length $L^*$, e.g. that there is no closed (retracing) walk on $G$ which can be done in less than $L^*$ steps.
\subsection{Cycles and cuts}
As we will see in the next sections, cycles can be reduced via Cauchy-Schwartz inequality (the cuts mentioned below) to the analysis of a single cycle. For this reason, it is worth first understanding what happens when the interaction graph has a single cycle.
\subsubsection{One-cycle expansion}
Let us assume, to begin with, that we have a forest $\mathcal A$ and that there is, without loss of generality, a unique path augmentation between $i$ and $j$ with a cycle of length $L^\prime$. Such expansion is exact on graphs which have only one cycle.
Let us call $\mathcal P_{ij}\subset \mathcal F$ and $\mathcal A\setminus \mathcal P_{ij}$. For $\Delta \mathcal Z=(Z-Z_{\mathcal F})|_{1 cycle}$, we have
 \begin{eqnarray}
   |\Delta Z| &=&\sum_{\mathcal A} \int d\mu(\sigma) \prod_{(rt)\in \mathcal A\setminus \mathcal P_{ij}}  (e^{\beta A_{rt} H_{rt}}-1) \nonumber \\
    &\cdot&\int_0^1\cdots \int_0^1 \prod_{(rt)\in  \mathcal P_{ij}} du_{rt}\ (\beta A_{rt} H_{rt}) e^{\sum_{e_{kt} \in \mathcal P_{ij}} A_{rt} H_{rt} u_{rt} } \nonumber \\
    & \cdot& (e^{A_{ij} H_{ij} \text{min}(u_{rt}\in \mathcal P_{ij}) }-1).\nonumber 
\end{eqnarray}
 
We note that the introduction of any extra edge on the interaction graph, which connects two vertices $i,j\in \mathcal F$ in the same tree, by construction, introduces a cycle of length $L=|P_{ij}|+1$. 
 
 At this point we introduce in this simpler case the techniques to upper bound $|\Delta Z|$. We focus first on the 1-cycle expansion, and extend these results to many (non necessarily disconnected) cycles in the following.
 In particular, we are interested in obtaining upper bounds of the form
 \begin{eqnarray}
     |Z-Z_{\mathcal F}|\leq \mathcal S(L^*)
 \end{eqnarray}
 where $L^*$ is the minimum length of a cycle in the interaction graph, and
 where $\mathcal S(L^*)$ is a function which goes to zero for $L^*\rightarrow \infty$.

 In eqn. (\ref{eq:forestred}) there are various terms that appear and that need to be taken care of if we want to introduce a proper bound. First we focus on
 \begin{eqnarray}
     & &\int_0^1\cdots \int_0^1 \prod_{(rt)\in  \mathcal P_{ij}} du_{rt} \beta A_{rt} H_{rt} e^{\sum_{e_{kt} \in \mathcal P_{ij}} A_{rt} H_{rt} u_{rt}} 
     %\nonumber \\
     %& &\ \ \ \ \ \ \ \ \ \ \ \ \cdot  
     e^{ A_{ij} H_{ij} \text{min}(u_{rt}\in \mathcal P_{ij}) },
 \end{eqnarray}
 which the expression typically appearing if a cycle can be formed via the tree expansion.
 The technique we employ to separate the tree expansion from the cycles, is to surgically remove the cycle, associated with the ``min" term in the exponential. An upper bound can be in fact obtained via the Cauchy-Schwartz inequality for integrals, as we now explain. First, let us note that 
 \begin{eqnarray}
     \int_0^1\cdots \int_0^1 d u_1 \cdots du_k (\cdot) =\langle \cdot \rangle
 \end{eqnarray}
 is effectively a measure, and $\langle 1\rangle=1$. Let us call $\int Du^k$ such measure. Then, we observe that given two real functions $f(\vec u)$ and $g(\vec u)$, we can define
 \begin{eqnarray}
     \langle f(\vec u),g(\vec u)\rangle&=&\int Du^k f(\vec u) g(\vec u) \nonumber \\
     &=&\int_0^1 du_1\cdots du_k f(u_1,\cdots,u_k)g(u_1,\cdots,u_k)
     \nonumber 
 \end{eqnarray}
 
 \begin{figure}
    \centering
    \includegraphics[scale=1.4]{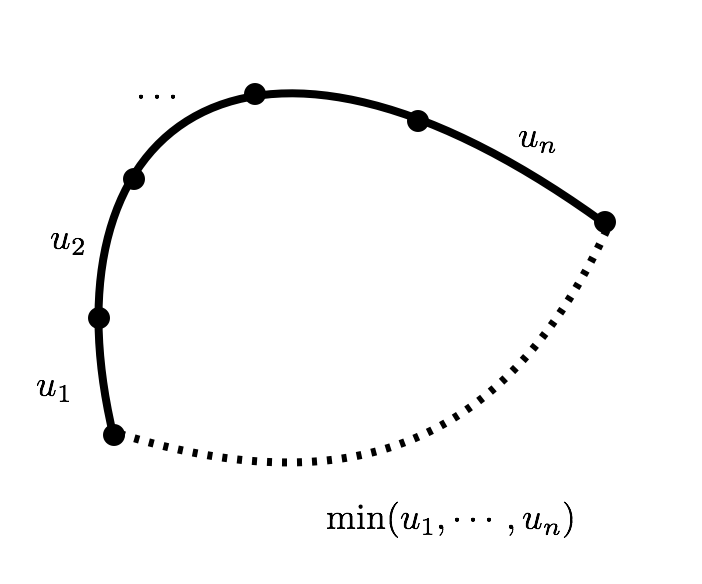}
    \caption{This is an example of a cycle and the closure variable involved.}
    \label{fig:singleloop}
\end{figure}
 
 The basic property of integrals that we will use is
 \begin{eqnarray}
     |\int_0^1 D u^k f(\vec u) g(u)|\leq \text{sup}_{\vec u} |f(u)|\cdot | \int_0^1\cdots \int_0^1 D u^k  g(u)|.\nonumber 
 \end{eqnarray}
 
Consider then
 \begin{eqnarray}
     \int_0^1&\cdots& \int_0^1 \prod_{(rt)\in  \mathcal P_{ij}} du_{rt} \beta A_{rt} H_{rt} e^{\beta \sum_{e_{kt} \in \mathcal P_{ij}} A_{rt} H_{rt} u_{rt}}   \nonumber \\
     &\cdot& (e^{\beta A_{ij} H_{ij} \text{min}(u_{rt}\in \mathcal P_{ij}) })-1)\nonumber \\
     &=&\langle f_{\mathcal F}(\vec u),g_{\mathcal F}(\vec u)\rangle,
 \end{eqnarray}
 where we choose
 \begin{eqnarray}
    f_{\mathcal F}(\vec u)
    &=&  \prod_{(rt)\in  \mathcal P_{ij}}  \beta A_{rt} H_{rt} e^{\sum_{e_{kt} \in \mathcal P_{ij}} A_{rt} H_{rt} u_{rt}} \nonumber \\
    g_{\mathcal F}(\vec u)&=&e^{\beta A_{ij} H_{ij} \text{min}(u_{rt}\in \mathcal P_{ij})}-1
 \end{eqnarray}
 The situation is similar to the one of Fig. \ref{fig:singleloop}.
 Then, we have
 \begin{eqnarray}
     \langle f_{\mathcal F},g_{\mathcal F}\rangle&\leq& sup_{\vec u}|\prod_{(rt)\in  \mathcal P_{ij}}  \beta A_{rt} H_{rt} e^{\sum_{e_{kt} \in \mathcal P_{ij}} A_{rt} H_{rt} u_{rt}}|\nonumber \\
     & &\cdot |\int_0^1 D u^k e^{\beta A_{ij} H_{ij} \text{min}(u \in \mathcal P_{ij})}-1| \nonumber \\
     \label{eq:ineq1} 
 \end{eqnarray}

We will focus on the integral above next. Let us now see how to bound these integrals when we have multi-cuts.

\begin{figure}
     \centering
     \includegraphics[scale=1.4]{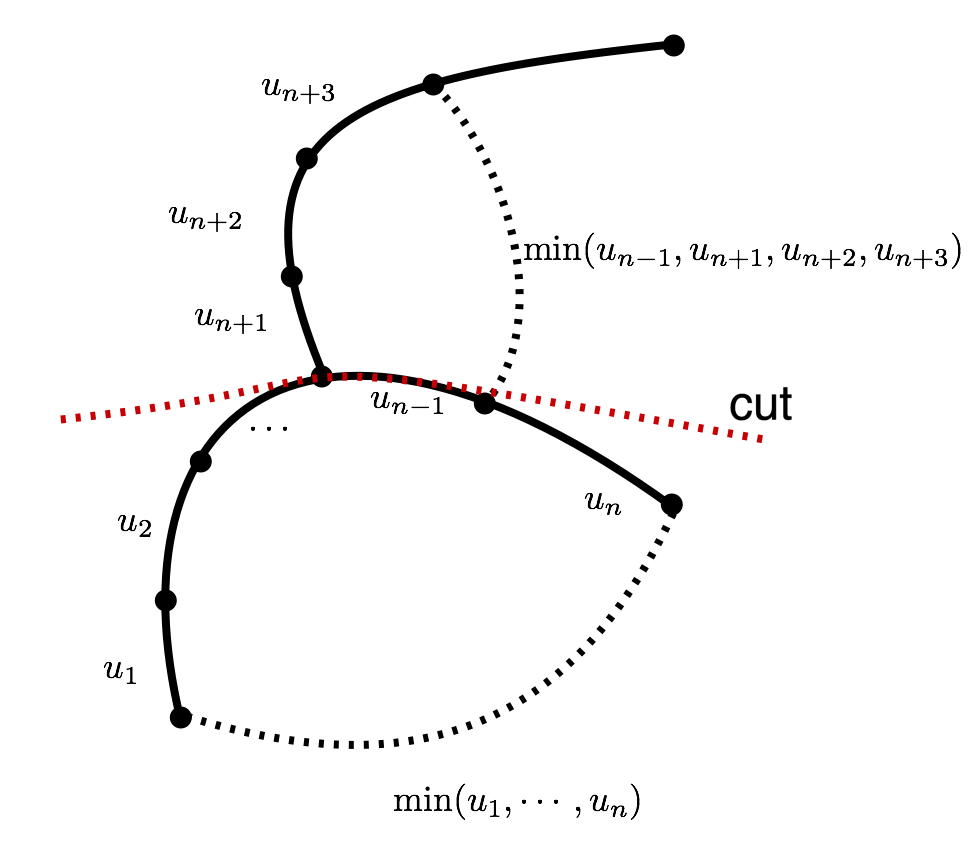}
     \caption{Example of a cut in the case of a multi-cycle forest element.}
     \label{fig:multiloops}
 \end{figure}

 It is easy to show that for semi-positive real functions, $\langle f,g\rangle$ of eqn. (\ref{eq:ineq1} ) is a scalar product. We have in fact
 \begin{itemize}
     \item $\langle f(\vec u),f(\vec u)\rangle=\int Du^k f^2(\vec u)\geq0$, and thus $\langle f,f\rangle=0 \implies f=0$.
     \item The scalar product is clearly bilinear, as $\langle a f,g\rangle=a \langle f,g\rangle$ and  $\langle  f,a g\rangle=a \langle f,g\rangle$.
     \item The property of conjugate definiteness is trivial for real functions.
 \end{itemize}
 This implies that we can use the Cauchy-Schwartz Inequality (CSI), and we have in fact
 \begin{eqnarray}
     \langle f,g\rangle\leq \sqrt{\langle f^2\rangle\langle g^2\rangle}.
 \end{eqnarray}
 Let us now introduce the notion of cuts, which is based on the CSI and have been introduced in \cite{magnen}. Consider the situation in which there are multiple paths in the forest expansion, a situation with two cycles as follows
  \begin{eqnarray}
     \langle f_{\mathcal F},g_{\mathcal F}\rangle&\leq& \text{sup}_{\vec u}|\prod_{(rt)\in  \mathcal P_{ij}}   \beta A_{rt} H_{rt} e^{\beta\sum_{e_{kt} \in \mathcal P_{ij}} A_{rt} H_{rt} u_{rt}}|\nonumber \\
     &\cdot& |\int_0^1\cdots \int_0^1 du_1\cdots u_k e^{\beta A_{ij} H_{ij} \text{min}(u \in \mathcal P_{ij})}-1| \nonumber \\
 \end{eqnarray}
 we will see later that the case with $N$ cycles can be reduced to the case with $2$ cycles via multiple cuts.  The situation we are interested in particular is like the one in Fig. \ref{fig:multiloops}, in which two cycles have a number of variables in common.
 
 In the following we will assume that interactions are bounded, e.g. that
 \begin{eqnarray}
     \bar H=\text{sup}_{\vec u}|\prod_{(rt)\in  \mathcal P_{ij}}  \beta A_{rt} H_{rt} e^{\sum_{e_{kt} \in \mathcal P_{ij}} A_{rt} H_{rt} u_{rt}}|< \infty.\nonumber
 \end{eqnarray}
 In particular, what is important to note is that
 \begin{eqnarray}
     \lim_{L\rightarrow\infty} \frac{\bar H}{
     L}=0.
 \end{eqnarray}
 Now note that we can upper bound $\bar H$ as
 \begin{eqnarray}
     \bar H\leq \beta^{L-1}  (\text{sup} |H_{ij}|)^{L-1} e^{(L-1) \beta \text{sup} |H_{ij}|}.
 \end{eqnarray}
 We see that we can write 
 \begin{eqnarray}
     \beta^{L-1}  (\text{sup} |H_{ij}|)^{L-1} e^{(L-1) \beta \text{sup} |H_{ij}|}= (x e^x)^{L-1}
 \end{eqnarray}
 where $x=\beta \text{sup} |H_{ij}|$. It follows that there is a critical value for $x$ for which the approximation applies. This is given by
 \begin{eqnarray}
     x_c=W(1)\approx 0.567134[..],
 \end{eqnarray}
 where $W(x)$ is the Lambert W-function. It follows that for 
 \begin{eqnarray}
     T^*\geq \kappa W(1) \text{sup}{|H_{ij}|}
 \end{eqnarray}
 such approximation applies. In the following we set $\kappa=1$ and assume that such requirement holds. The value of $T^*$ thus marks point at which the system goes from a low temperature phase into a forest phase (in the limit $L^*\rightarrow \infty$).

 \subsubsection{Cycles, Cuts and Cauchy-Schwartz inequality}

We now discuss many cycles in the same forest. Let us first understand how cycles enter into the calculation.
The goal will be to bound the contribution of the cycles into the partition function in the size of the minimum length $L$. Let us first clarify the various cases that can occur, and we will analyze them separately. The goal is to reduce the upper bound to a sequence of 1-cycles which we studied in the previous section.

Situation \textit{(a)} is the one in which the cycles are on two different paths, e.g. $P_{ij}\cap P_{i^\prime j^\prime}=\{\emptyset\}$.
In this case we can split the integrals in two, as from the point of view of the variables $u_{ij}$ the two cycles are not overlapping on any variable $u_{ij}$. Thus $\text{min}(u_{ij}\in \mathcal P_{ij})$ and $\text{min}(u_{ij}\in \mathcal P_{i^\prime j^\prime })$ act on different variables. The situation \textit{(b)} is that the cycles overlap, meaning $P_{ij}\cap P_{i^\prime j^\prime}\neq \{\emptyset\}$.  An extreme situation of this is one in which \textit{(b1)} $P_{ij}\cap P_{i^\prime j^\prime}=P_{ij}$ or 
\textit{(b2)} $P_{ij}\cap P_{i^\prime j^\prime}=P_{i^\prime j^\prime}$. In this situation, one cycle is contained into the other. We can write  $ \text{min}(u_{ij} \in \mathcal P_{i^\prime j^\prime })= \text{min}(\text{min}(u \in \mathcal P_{i j }), u\in \{\mathcal P_{i^\prime j^\prime }\setminus \mathcal P_{ij}\})$ or
$ \text{min}(u_{ij} \in \mathcal P_{i j })= \text{min}(\text{min}(u \in \mathcal P_{i^\prime j^\prime }), u\in \{\mathcal P_{ij}\setminus \mathcal P_{i^\prime j^\prime }\}) $. The less extreme situation is the one in which the two augmented paths contain only a subset of the variables. Both situation can be overcome with a cut, which is an application of the Cauchy-Schwartz inequality in order to separate two cycles, and introduced first in \cite{Gurau2,magnen}. The particular technique we employ in this paper has been introduced in \cite{Caravelli}. Specifically, we have the following general expression and upper bound for the situation of type (b):

\begin{widetext}
\begin{eqnarray}
    |\int_0^1\cdots \int_0^1 \prod_{(rt)\in  \{\mathcal P_{ij}\cup\mathcal P_{i^\prime j^\prime}\}}& & du_{rt} \beta A_{rt} H_{rt} e^{\sum_{e_{kt} \in \{\mathcal P_{ij}\cup \in \mathcal P_{i^\prime j^\prime}\}} \beta A_{rt} H_{rt} u_{rt}} \nonumber \\
    &\cdot& (e^{\beta A_{ij} H_{ij} \text{min}(u_{rt}\in \mathcal P_{ij})}e^{\beta A_{i^\prime j^\prime } H_{i^\prime j^\prime} \text{min}(u_{rt}\in \mathcal P_{i^\prime j^\prime}) }-1)|\nonumber \\
    &\leq& \text{sup}_{\vec u}|\beta A_{rt} H_{rt} e^{\sum_{e_{kt} \in \{\mathcal P_{ij}\cup \in \mathcal P_{i^\prime j^\prime}\}} \beta A_{rt} H_{rt} u_{rt}} | \nonumber \\
    & & \ \ \ \ \ \ \cdot |\int_0^1 du^k(e^{\beta A_{ij} H_{ij} \text{min}(u_{rt}\in \mathcal P_{ij})}e^{\beta A_{i^\prime j^\prime } H_{i^\prime j^\prime} \text{min}(u_{rt}\in \mathcal P_{i^\prime j^\prime}) }-1)|.
\end{eqnarray} 
\end{widetext}
In order to clarify the sequence of cuts in what follows, we introduce the following notation. Let us now call the average for $K$ cycles as
\begin{eqnarray}
    Q_K^T(L_1,\cdots,L_K)=\int_0^1 du^k\  \prod_{k=1}^Ke^{T \beta A_{i_k j_k} H_{i_kj_k} \text{min}(u_{rt}\in \mathcal P_{i_kj_k})  }\nonumber \\
\end{eqnarray}
where we identify with $L_1$ the length of the cycle, e.g. we have in each subcycle the variables $u_{i_1}\cdots u_{i_{L-1}}$.
We focus on 
\begin{eqnarray}
    Q_2^1(L_1,L_2)&=&\int_0^1 du^k\ e^{\beta A_{ij} H_{ij} \text{min}(u_{rt}\in \mathcal P_{ij})} \nonumber\\
    & & \ \ \ \ \ \cdot e^{\beta A_{i^\prime j^\prime } H_{i^\prime j^\prime} \text{min}(u_{rt}\in \mathcal P_{i^\prime j^\prime}) }
\end{eqnarray}
 We can use the Cauchy-Schwartz inequality as follows.
 We define
 \begin{eqnarray}
     f_{\mathcal F}(\vec u)&=&e^{\beta A_{ij} H_{ij} \text{min}(u_{rt}\in \mathcal P_{ij})}\nonumber \\
     g_{\mathcal F}(\vec u)&=&e^{\beta A_{i^\prime j^\prime } H_{i^\prime j^\prime} \text{min}(u_{rt}\in \mathcal P_{i^\prime j^\prime}) }
 \end{eqnarray}
 Then, we can apply the Cauchy-Schwartz inequality and obtain
 \begin{eqnarray}
     \langle f_{\mathcal F}(\vec u),g_{\mathcal F}(\vec u)\rangle\leq \sqrt{ \langle f_{\mathcal F}(\vec u),f_{\mathcal F}(\vec u)\rangle\langle g_{\mathcal F}(\vec u),g_{\mathcal F}(\vec u)\rangle} \nonumber 
 \end{eqnarray}
Evaluating these integrals $\langle g_{\mathcal F}(\vec u),g_{\mathcal F}(\vec u)\rangle$ and $\langle f_{\mathcal F}(\vec u),f_{\mathcal F}(\vec u)\rangle$ independently, we see that we have effectively decoupled the cycles, but doubled the effective strength of the interaction $A_{ij}$ and $A_{i^\prime j^\prime}$. We can apply for each single cycle the inequality of eqn. (\ref{eq:ineq1}) twice, which gives the same common factor squared, but with an overall square root.
Thus we reduce for $2$-cycles to the computation of 
\begin{eqnarray}
   |Q^1_2(L_1,L_2)|&\leq&  \sqrt{ \int_0^1 Du^k e^{2\beta A_{ij} H_{ij} \text{min}(u_{rt}\in \mathcal P_{ij})}  } \nonumber \\
   &\cdot& \sqrt{\int_0^1 D\tilde u^k e^{2\beta A_{i^\prime j^\prime } H_{i^\prime j^\prime } \text{min}(\tilde u \in \mathcal P_{i^\prime j^\prime })}}\nonumber \\
   &=&\sqrt{Q_1^2(L_1)Q_1^2(L_2)}
\end{eqnarray}
which is twice the calculation for a single cycle.

At this point we can proceed in the calculation showing that multiple overlapping cycles can be reduced to two overlapping cycles via CSI, as done in \cite{Caravelli}. We thus consider the case in which there are multiple cycles in the forest expansion. It is not hard to see that we can write, via a sequence of cuts,
\begin{eqnarray}
    |Q_K^1(L_i)|&\leq& \sqrt{Q_1^2(L_1)Q_{K-1}^2(L_2,\cdots,L_K)} \nonumber \\
    &\leq& \sqrt{Q_1^2(L_1)\sqrt{Q_1^{4}(L_2)Q_{K-2}^4(L_3,\cdots,L_K)}}\nonumber \\
    &\vdots&\\
    &\leq &\prod_{k=1}^K (Q^{2k}_1(L_k))^{\frac{1}{2k}}.
\end{eqnarray}
Effectively, we  have now reduced the multi-cycle case, via a sequence of cuts, to a 1-cycle calculation which we know how to resolve as we will see shortly. The key issue at this point is how many cycles we should expect, and whether the limit $L\rightarrow \infty$ does affect the upper bound. For this a little care is necessary.  
\subsubsection{Cycle Sparse versus Cycle Dense graphs}
\begin{figure}
    \centering
    \includegraphics[scale=1.2]{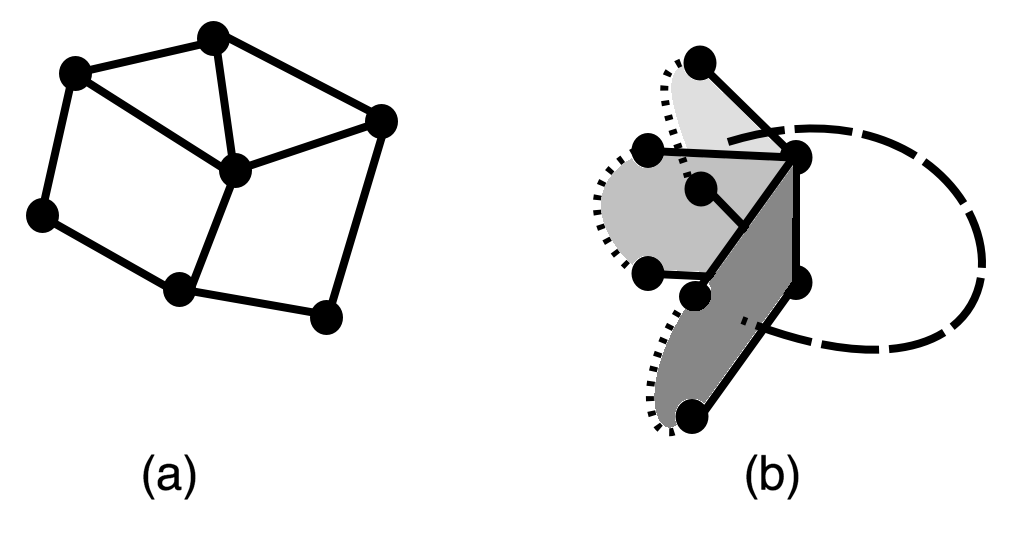}
    \caption{Two examples for cycle dense and cycle sparse graphs. In (a) we see that each cycle is only finitely connected to many other cycles via connected links. In (b) we have that each cycle is connected to all the other $K$ cycles via a single edge, independently from the length of the cycle. In this case, if these cycles are of length $L_1\cdots L_K$, then the total number of edges of the graph is $E=\sum_{i=1}^K L_i-K+1$, and if the cycles are all equal, $E=K(L-1)+1$. Since these are all cycles, we have also $N=K(L-2)+2$. Thus $E-N=K-1$ and $K=\frac{N-2}{L-2}$, which corresponds to $\eta=1$ in eqn. (\ref{eq:eta}).}
    \label{fig:loopdense}
\end{figure}

As we have seen, multiple adjacent cycles can be reduced to a single cycle. This said, two different analysis apply, depending on the scaling of the number of adjacent cycles, $K$ and most importantly how number scales with $L$. We wish to show here that the upper bound can depend on this scaling, without however affecting the final result, e.g. the fact that the forest expansion applies.
The first is the case for which $K$ is finite, and $\frac{K}{L}\rightarrow 0$. The second is when $K\rightarrow \infty$ and $L\rightarrow \infty$. It is easy to see that according to the bounds from the previous section, a graph is dense or sparse according to the limit of $K$ versus $L$. The value of $K$ can be defined via the dual graph. The dual graph $\mathcal L$ is the graph by which we replace every cycle by a vertex, and connect these vertices if two cycles are adjacent (e.g. they share an edge). Examples of cycle dense and sparse graphs are shown in Fig. \ref{fig:loopdense}.
\ \\\ \\
\textbf{Definition 1} Let $K$ be the maximum degree of the dual graph of $\mathcal G$.
We say that a graph $G$ is \textit{cycle sparse} if for $N\rightarrow \infty$, $\frac{K}{L}\rightarrow 0$. On the other hand if $\lim_{N\rightarrow \infty} \frac{K}{L}\rightarrow c$ where $c$ is constant and nonzero, we say that the graph is \textit{cycle dense}.\ \\\ \\

The definitions are provided below.

%What we wish to show is that what discriminates the quality of the forest expansion is exactly the sparse graph regime.
As we see above, the key issue is how the upper bound behaves when we perform the limit $L^*\rightarrow \infty$. As we will see below, if $K$ is finite (only a finite set of loops are adjacent), then the limit is well defined and we can perform the expansion without further analysis. If however the graph is dense, then we need to take extra care in performing the limit. However, as we will see, we can still make sense of the forest expansion.

\textit{Cycle Sparse graphs: Case $K$ finite, $L$ large.}

First let us assume $K$ finite. 
Let us anticipate that for cycle sparse graphs we have
\begin{eqnarray}
    Q_1^{2k}(\beta A_{ij} H_{ij},L_1)\approx 1+\frac{ 2k \beta A_{ij} H_{ij}}{L_1}+O(\frac{1}{L_1^2}).
\end{eqnarray}

 If we use this formula, we have 
\begin{eqnarray}
    |Q_K^1(L_1,\cdots,L_K)|\leq \prod_{k=1}^K (1+\frac{2k \beta A_{ij} H_{ij}}{L_k})^{\frac{1}{2k}}.
\end{eqnarray}

If we now assume that $L_k\gg \beta H_{ij}$, we can write
\begin{eqnarray}
\prod_{k=1}^K (1+\frac{2k \beta A_{ij} H_{ij}}{L_k})^{\frac{1}{2k}}\approx     \prod_{k=1}^K (1+\frac{ \beta A_{i_k j_{k}}H_{i_k j_{k}}}{L_k})
\end{eqnarray}
If we want to keep only the terms of order up to $\frac{1}{L}$, then we have 
\begin{eqnarray}
    \prod_{k=1}^K (1+\frac{A_{i_k j_{k}}H_{i_k j_{k}}}{L_k})^{\frac{1}{2k}}\approx 1+\sum_{j=1}^K \frac{\beta A_{i_k j_{k}}H_{i_k j_{k}}}{L_j}
\end{eqnarray}
from which we obtain a loose upper bound of the type
\begin{eqnarray}
    |Q_K^1(L_1,\cdots ,L_K)| \leq  1+  \beta\bar H\frac{K}{L^*}
\end{eqnarray}
where $L^*$ is the minimum cycle size in the interaction graph, and $\bar H= \text{sup}_{k}| H_{i_k j_k}|$.
The result above applies only if $K$ is finite, e.g. if $K$ does not scale with $L$ (or $N$).

The bound we just obtained relied on a specific identity:
\begin{eqnarray}
    Q_1(A)&=&\int_0^1\cdots \int_0^1 e^{ A \text{min}(u_1,\cdots,u_L)} du_1\cdots u_L \nonumber \\
    &=&\frac{\Gamma(L+1)}{\Gamma(L)}e^A A^{-L} (\Gamma (L)-\Gamma (L,A)) \nonumber \\
    &\approx& 1+\frac{A}{L+1}+O(\frac{1}{L^2})
    \label{eq:q1}
\end{eqnarray}
which we now prove.
%\textit{Proof}
First we note that by the symmetry of the integral, we can write
\begin{eqnarray}
   S_L(A)&=& \int_0^1\cdots \int_0^1 e^{ A \text{min}(u_1,\cdots,u_L)} du_1\cdots u_L \nonumber \\
   &=&  \Gamma(L+1) \int_0^1 du_1 \int_0^{u_1}\cdots \int_0^{u_{L-1}} du_L e^{ A u_L} \nonumber.
\end{eqnarray}

We note that the formula can be obtained via the Cauchy recursive integration of the integrals. Given a function $f(x)$, we define
\begin{eqnarray}
    f^{(-n)}(x)=\int_a^x \int_a^{u_1} \cdots \int_a^{u_{L-1}} f(u_L) du_L\cdots du_1.
\end{eqnarray}
Cauchy formula shows that we can evaluate the integral via the calculation of a reduced one given by
\begin{eqnarray}
    f^{(-n)}(x)=\frac{1}{\Gamma(L)}\int_a^x (x-t)^{L-1} f(t)dt.
\end{eqnarray}
We thus have
\begin{eqnarray}
    S_L(A)&=&\int_0^1\cdots \int_0^1 e^{ A \text{min}(u_1,\cdots,u_L)} du_1\cdots u_L \nonumber \\
    &=&\frac{\Gamma(L+1)}{\Gamma(L)}\int_0^1 (1-t)^{L-1} e^{At} dt \nonumber \\
    &=&\frac{\Gamma(L+1)}{\Gamma(L)}e^A A^{-L} (\Gamma (L)-\Gamma (L,A))
\end{eqnarray}
where $\Gamma(L,A)=\int_A^\infty t^{L-1} e^{-t} dt$.
However, what we care about is the leading order in $L$ for $L\rightarrow \infty$. We have
\begin{eqnarray}
    \frac{\Gamma(L+1)}{\Gamma(L)}e^A A^{-L}& &(\Gamma (L)-\Gamma (L,A)) \nonumber \\
    &\sim& 1+\frac{A}{L+1}+\frac{A^2}{(L+2)(L+1)} \nonumber \\
    &+&\frac{A^3}{(L+3)(L+2)(L+1)}+\cdots.\nonumber
\end{eqnarray}
which concludes the proof.%$\square$

\textit{Cycle Dense graphs: Case $K$ and $L$ infinite.}
Let us now consider the case in which now $K$ is unbounded, but scales with the size of the loop in a particular manner.
First, note that in the worst case scenario, these cycles are connected via a $1$-edge each. Thus, there is a number of $u$ variables involved given by $N_u=L_1+\cdots+L_K-K$. Since the number of $u$ variables has to be $N_u\leq \frac{N(N-1)}{2}$, we must have that for the minimum length of the cycle we have at most $K(L^*-1)\leq \frac{N(N-1)}{2}$. In general we can parametrize the growth of $K$ with $N$ via the scaling
\begin{eqnarray}
    K\sim \frac{N^{\eta}}{  L^*},
    \label{eq:eta}
\end{eqnarray}
with $0\leq \eta\leq 2$ for $K\rightarrow \infty$, $N\rightarrow \infty$ and $N\rightarrow \infty$. On the other hand, we have the general bound 
\begin{eqnarray}
    3\leq L^*\leq N-1.
\end{eqnarray}
Since we are interested in $L^{*}\rightarrow \infty$, necessarily we must have $L^*= N^\gamma$ for $c$ and $\gamma$ constants. Analogously, assume $K\sim  N^{\gamma^\prime}$. Where we must have via a scaling argument that
\begin{eqnarray}
    \gamma = \eta- \gamma^\prime.
\end{eqnarray}
and $L^*\sim K^{\frac{\gamma}{\gamma^\prime}} $.
It follows that we can perform the following upper bound
\begin{eqnarray}
    \prod_{k=1}^K (1+\frac{A_{i_k j_{k}}H_{i_k j_{k}}}{L_k})^{\frac{1}{2k}}&\leq& \prod_{k=1}^K (1+\frac{2k \beta A_{ij} H_{ij}}{L^*})^{\frac{1}{2k}}\nonumber \\
    &\leq&  \lim_{K\rightarrow \infty} (1+\rho \frac{\beta \bar H}{K^{\frac{\gamma}{\gamma^\prime}}})^K,
\end{eqnarray}
where $\rho$ depends on the scaling relationships for $L^*$ and $K$ in terms of $N$.
Because of the relationship above, we have $\frac{\gamma}{\gamma^\prime}=\frac{\eta}{\gamma^\prime}-1$. Now note that if $\frac{\eta}{\gamma^\prime}>2$, then such limit is one, while if $\frac{\eta}{\gamma^\prime}=2$ the limit is $e^{\rho \beta H}$.
We can now take the limit $K\rightarrow \infty$, from which we obtain
\begin{eqnarray}
    \prod_{k=1}^\infty  (1+\frac{A_{i_k j_{k}}H_{i_k j_{k}}}{L_k})^{\frac{1}{2k}}\leq \lim_{K\rightarrow \infty}\prod_{k=1}^K (1+\frac{\rho \beta  \bar H}{K})=  e^{\rho  \beta \bar H }.\nonumber
\end{eqnarray}

The argument above shows that, since the cycle contributions are finite (if H is bounded) even in the case of dense cycles the forest approximation  works. In fact, $e^{\rho \beta \bar H}-1\approx \rho \beta \bar H$ is finite in this case, and the two limits are parametrized by $\rho\geq 0$.

\subsection{Sum over $\mathcal A_K$ and Forest formula corrections}
In order to simplify the calculations that follow, we note that we can write the expansion over $\mathcal A$ in terms of a sum over all graphs which contain $K$ cycles. This is $\sum_{\mathcal A}=\sum_{K=1} \sum_{\mathcal A_{\mathcal K}}$ without loss of generality, where  $\mathcal A_{\mathcal K}$ a graph which contains $K$ path augmentations. Each element in $\mathcal A_{\mathcal K}$ can be written as a graph sum such that there are $K$ path augmentation, which is essentially as fixing the cycles in the graph. The sum is on every forest of the remainder of the graph which contains that specific combination of path augmentations.

If we put together the bounds from the previous section, we obtain
\begin{widetext}
\begin{eqnarray}
    |Z-Z_{\mathcal F}|\leq \int d\mu(\sigma)\Big(\sum_{K}K g(L)\sum_{\mathcal A_K}  \prod_{(rt)\in  {\mathcal A}_K\setminus \cup_{k} P_{i_k j_k}} (e^{\beta A_{rt} H_{rt}}-1)\cdot \prod_{k=1}^K   (\beta \bar H)^{|P_{i_k j_k}|-1} e^{\beta \bar H(|P_{i_k j_k}|-1)} \Big)\nonumber \\
    \label{eq:forestbound}
\end{eqnarray}
\end{widetext}
where we defined $A_{\mathcal K}$ as the set of trees which have $K$-path augmentations (e.g. $K$ cycles), and $|P_{i_k j_k}|$ is the length of the cycle associated with the cotree element $(i_k j_k)$. Also, we have introduced the function $g(L)$ which is $\frac{K \beta \bar H}{L}$ for cycle sparse graphs and $e^{\rho^\prime \bar H}$ for cycle dense graphs, where $K$ is the number of cycles and $\rho$ is a parameter which characterizes how dense the graph is.

The proof would be done at this stage if we knew that $\sum_{\mathcal A_K}$ involves only a small number of elements. Unfortunately, Cayley's theorem shows that the number of such elements is humongously large, e.g.  $\sum_{\mathcal A_K} \sim N^N$ elements. Thus, any attempt at performing an upper bound  would be hopeless at this stage, even if this multiplies something which is exponentially small.
We thus need to discuss how to perform the sum in order to show that this is less than exponentially large in $L$. This issue requires some analysis on the meaning of the formula of eqn. (\ref{eq:forestbound}) in terms of cycles and forest, and will be discussed in the next section in more detail.

In particular, we will need to show that for $T\gg T^*$, then
\begin{eqnarray}
    \sum_{\mathcal A_K}& &  \prod_{(rt)\in \mathcal A_K\setminus \cup_{k} P_{i_k j_k}} (e^{\beta A_{rt} H_{rt}}-1) \nonumber \\
    &=& e^{-\beta H}|_{\mathcal T_{\mathcal A_K}}+ O\big( (\beta \bar H)^{\sum_k |P_{i_k j_k}\in \mathcal A_{\mathcal K}|} \big)
    \label{eq:finalproof}
\end{eqnarray}
where $\mathcal T_{\mathcal A_K}$ is the tree obtained from removing the $K$ cycles associated to the path augmentations from the graphs.

At this stage, all we can say is that if we can prove that $\sum_{\mathcal A_K}  \prod_{(rt)\in \mathcal F\setminus \cup_{k} P_{i_k j_k}} (e^{\beta A_{rt} H_{rt}}-1)$ is bounded, then
\begin{eqnarray}
    Z(\beta)=Z_{\mathcal F}+\begin{cases}
    O(q^{L-1} e^{\rho \beta \bar H}) & \text{cycle dense graphs} \\
   O( \frac{q^{L-1}}{L}) & \text{cycle sparse graphs}
    \end{cases}\nonumber
\end{eqnarray}
In either case, we see that what determines a high and a low temperature is given by the value of $q$.
The value of $q$ is temperature dependent, if we define $\bar H=\text{sup}_{ij} |H_{ij}|$, then
\begin{eqnarray}
    q=\beta \bar H e^{\beta H}.
\end{eqnarray}
We thus have that $q<1$ for $\beta< \beta^*$, with $\beta^*=(W(1) \bar H)^{-1}$.$\square$

As a result, because of the statement above, we know that for $\beta<\beta^*$ we have
\begin{eqnarray}
    Z\approx Z_{\mathcal F}.
\end{eqnarray}
where $Z_{\mathcal F}$ is a polynomial expansion in the graph edge variables $x_{ij}=e^{\beta A_{ij} H_{ij}}-1$. The interpretation of the result is the one we presented in the results section at the beginning of the paper.
In the next section we clarify the role of spanning trees in the forest expansion, and provide the proofs for the second and third results claimed in the introduction.

\section{Expression for $Z_{\mathcal F}$ in terms of trees and cycles}
\subsection{Reference tree expansion}
We have seen that for sparse graphs with large cycles the probability distribution can be approximated by a forest expansion. 
In principle, every single term of this expansion can evaluated in polynomial time, being the expansion a product of variables $x_{ij}=e^{\beta A_{ij} H_{ij}}-1$ on a forest subgraph of $\mathcal G$. However, we wonder whether such expansion can be rewritten in terms of trees in a more compact form. Moreover, it is very well known that the number of forests on a graph on $N$ vertices grows as $(N-1)^{N-1}$, which is Cayley's formula. It is thus interesting to obtain closed form expressions for the forest expansion. In the following we work at the level of the probability distribution, but all results about the partition function can be obtained by integrating over the measure $\int D\mu(\sigma)$.

More specifically, we are interested in finding an expression for the partition function $Z_{\mathcal F}$. This is motivated by the fact that if we can write $Z_{\mathcal  F}$ in terms of trees, then its computation can be done efficiently via belief propagation, if not exactly, at least numerically. Note that graphs are vector spaces, and thus collections of graphs can be ``orthogonalized".

Now note that since in $Z_{\mathcal F}$ we have integrated out all $u_{ij}$, the result is only dependent on polynomials in $x_{ij}=e^{\beta A_{ij} H_{ij}}-1$. Given a graph $\mathcal G$ we define the graph polynomial
\begin{eqnarray}
    P_{\mathcal G}(x)=\prod_{(ij)}(1+x_{ij})
\end{eqnarray}
Let us introduce the following notation: given a graph $\mathcal G$, $\mathcal F(\mathcal G)$ is the sum over the forests, e.g. any connected or disconnected tree subgraph of $\mathcal G$.
We have that $Z_{\mathcal F}=\sum_{G\in \mathcal F(\mathcal G)} P_{ \mathcal G}(\vec x)$.
Now we note that one can write
\begin{eqnarray}
    e^{-\beta H}&=&e^{\beta \sum_{ij} A_{ij }H_{ij}}=\prod_{(ij)\in E(\mathcal G)} (1+(e^{\beta A_{ij} H_{ij}}-1)) \nonumber \\
    &=&\prod_{(ij)\in \mathcal G} (1+x_{ij})=\sum_{\mathcal G^\prime \subset \mathcal G} \prod_{(ij) \in  E(\mathcal G^\prime)} x_{ij} .
\end{eqnarray}
Thus the probability distribution can be written as a sum over all possible subgraphs of $\mathcal G$.
  Clearly, if $\mathcal G$ contains cycles, the sum is over all possible subgraphs (with cycles or not).  However, if $\mathcal G$ is a tree  all possible subgraphs define a forest sum. 
Thus, we have
\begin{eqnarray}
    \prod_{(ij)\in \mathcal T}(1+x_{ij})=\sum_{\mathcal G^\prime \in \mathcal F(\mathcal T)} \prod_{(ij)\in \mathcal G^{\prime}}x_{ij}.\label{Eq:treesz}
\end{eqnarray}
This implies that for trees we have $\mathcal Z_{\mathcal F}=\sum_{\mathcal G^\prime \in \mathcal F(\mathcal T)} P_{\mathcal G^\prime}(x)$, and thus simply we have confirmed the BKAR formula for the case of a tree graph. We write such identification as 
    $\mathcal Z_{ {\mathcal F}{(\mathcal T)}}=P_{\mathcal F(\mathcal T)}(\{x\})=\mathcal Z_{\mathcal T}$,
e.g. $Z_{\mathcal T}$ is simply the partition function of the tree. Such expansion is commonly used to show that, for instance, belief propagation or cavity methods are exact on trees.

Clearly, the formula is not true for the case in which $\mathcal G$ contains cycles, as the expansion contains extra terms.
The message of the analysis above, albeit trivial, suggests that it is always possible to write $\mathcal Z_{\mathcal F}$ in terms of a reference tree (chosen arbitrarily), and analyse the leftover terms separately. 
Also, one comment that we can make immediately is that the $\mathcal Z_{\bar A}$ (e.g. the rest of the expansion), while being a sum over graphs without cycles, does contain information on the cycles as we saw. In fact, the forest will contain trees which surround a certain cycle without forming one. If we can write $\mathcal Z_{\mathcal F}$ in terms of trees. Clearly, such problem does not have naturally a unique solution, as given an arbitrary graph we have a number of spanning trees given by the Kirchhoff Matrix-Tree theorem. However, we wonder whether we can express $\mathcal Z_{\mathcal {\mathcal F}}$ in terms of a reference spanning tree $\mathcal T$.

 Let $\mathcal T$ be an arbitrary spanning tree $\mathcal T\subset \mathcal G$, and we define the co-tree as $\bar {\mathcal T}=\mathcal G\setminus \mathcal T$. Given an edge of the $e$ co-tree, we associate the cycle $C_e \in \mathcal G$ which contains $e$. We define the reduced cycle $\tilde C_e\setminus \{e\}$, which is the tree obtained from cycle by removing $e$. Clearly $\tilde C_e\subset \mathcal T$.
 It is known that $\mathcal C=\{C_{e}: e\in \bar {\mathcal T}\}$ forms a basis for all the cycles in the graph, and we will use this fact later. Also, we call $F_{e_1 \cdots e_k}(\mathcal G)$ a $k$-rooted forest of $\mathcal G$, e.g. all the forests which contain the edges $e_1,\cdots, e_k$. Given the above, we now introduce the following definition:\ \\\ \\
 
 \textit{Orthogonalization of rooted $k$-forests}. Let $e_1\cdots e_K$ be an ordered sequence of edges. We define the following sequence of reduced forests induced by such edge ordering:
 \begin{eqnarray}
    \tilde { \mathcal F}_{e_1 \cdots e_K}(\mathcal G)&=& { \mathcal F}_{e_1\cdots e_K}(\mathcal G) \nonumber \\
    \tilde { \mathcal F}_{e_1 \cdots e_{K-1}}(\mathcal G)&=& { \mathcal F}_{e_1 \cdots e_{K-1}}(\mathcal G)- \tilde { \mathcal F}_{e_1\cdots e_K}(\mathcal G)\nonumber \\
    \tilde { \mathcal F}_{e_1 \cdots e_{K-2}}(\mathcal G)&=& { \mathcal F}_{e_1 \cdots e_{K-2}}(\mathcal G)- \tilde { \mathcal F}_{e_1\cdots e_{K-1}}(\mathcal G)-\tilde { \mathcal F}_{e_1\cdots e_{K}}(\mathcal G)\nonumber \\
    &\vdots& \nonumber \\
    \tilde{\mathcal  F}_{e_1}(\mathcal G)&=&{\mathcal  F}_{e_1}(\mathcal G)-\tilde{\mathcal  F}_{e_1 e_2}(\mathcal G)-\cdots-\tilde { \mathcal F}_{e_1 \cdots e_K}(\mathcal G)
 \end{eqnarray}
  By construction, rooted forests are defined in such a way that they contain only the forests with the specified  labels of the relevant edges.
Thus, the orthogonalized forest are constructed such that $\tilde {\mathcal F}_{e_{i_1} \cdots e_{i_k}} (\mathcal G)$ contain only forests which contain the edges $e_{i_1} \cdots e_{i_k}$. It follows that:\ \\\ \\
 
  For an arbitrary graph $\mathcal G$ and a tree-cotree splitting $(\mathcal T,\bar {\mathcal T})$, with $|\bar {\mathcal T}|=K$
 \begin{eqnarray}
     \mathcal F(\mathcal G)&=&\mathcal F(\mathcal T)+\sum_{k=1}^K \sum_{\mathcal P_k(e_1,\cdots,e_K)}  \tilde{\mathcal F}_{e_{i_1} \cdots e_{i_k}}(\mathcal G) \nonumber \\
    &=&\mathcal F(\mathcal T) \nonumber \\
    &+&\sum_{k=1}^K \sum_{\mathcal P_k(e_1,\cdots,e_K)}  e_{i_1} \circ \cdots \circ e_{i_k}\circ \Delta F(e_{i_1},\cdots,e_{i_k}) \nonumber 
 \end{eqnarray}
 where the sum is over all $k$-partitions $\mathcal P_k(e_1,\cdots,e_K)$ of the edges $e_i\in \bar {\mathcal T}$, 
 and $\Delta F(e_{i_1},\cdots,e_{i_k})={\mathcal F}(\mathcal T)- Q_{e_{i_1}, \cdots, e_{i_k}}(\mathcal T)$ is the forest expansion without cycles.
 
 \textit{Proof}. 
 The proof of first equation follows from the fact that orthogonal rooted forests with different rooted edges necessarily contain different forest graphs, which is obtained from first statement. The second statement follows from the fact that all $k$-rooted forests are contained in $e_{i_1} \circ \cdots \circ e_{i_k}\circ \mathcal F(\mathcal T)$. However, from these one needs to subtract those graphs which generates cycles.  $\square$\ \\\ \\
 
 At this point, we can complete the proof of eqn. (\ref{eq:finalproof}) by inspecting what $\mathcal A_K$ represents. Let us assume that the total number of (fundamental) cycles is $P$. This is the sum over all possible graphs with $K\leq P$ cycles, and the cycles have been ``removed" and incorporated in the cycle upper bounds of the previous section. Thus, the sum $\mathcal A_K$ is essentially over all possible forests with the $K$ branches of the tree associated to the cycles cut off, but still (essentially) containing the remaining $P-K$ cycles, but failing to contain a path augmentation by construction. As such, similarly to what we had done before, we have in fact that for $\beta<\beta^*$, we have an expansion (or, alternatively, an upper bound on $|\mathcal Z_{\mathcal F}-\mathcal Z_{\mathcal T}|$) in terms of the reference tree expansion
 \begin{eqnarray}
    \sum_{\mathcal A_K}& &  \prod_{(rt)\in \mathcal A_K\setminus \cup_{k} P_{i_k j_k}} (e^{\beta A_{rt} H_{rt}}-1) %\nonumber \\     &\leq &  
     =e^{-\beta H}|_{\mathcal T_{\mathcal A_K}}+ O\big( \beta \bar H_{\mathcal G \setminus \mathcal T_{\mathcal A_K}} \big)
\end{eqnarray}
 where now $e^{-\beta H}|_{\mathcal T_{\mathcal A_K}}$ , which is the probability distribution of a reference tree, thus without cycles, and thus is of order zero with respect to $L^*$. Above, the notation $\bar H_{\mathcal G \setminus \mathcal T_{\mathcal A_K}}$ refers to the maximum coupling on the set of edges not in the reference tree.
 
 This concludes the proof of the convergence to the forest expansion of the previous section, as now in the limit $L^*\rightarrow \infty$, every single term multiplying $q^{(L^*-1)}$ is of order zero with respect to $L^*$.
 
 Since the rest does depend on the reference tree, minimizing the rest is a well defined optimization problem, as we discuss in the next section. It is interesting also to note that the correction is linear in $\beta$. Thus, as we discussed in the section at the beginning of the paper there is an intermediate ``forest phase" between $\bar H^{-1}\leq \beta \leq \beta^*$, and a cycle perturbative phase for $\beta>\bar H^{-1}$.
 
 \subsection{Optimal reference tree and cycle algebra}\label{sec:cyclealg}
 The last outstanding issue is how to choose the reference tree.
 Consider now orthogonal $1$-roots $\bar {\mathcal F}_{e_{i}}(e_{i}\circ \mathcal T)$.
 
 Let us  note that the previous construction gives a way of selecting the spanning tree $\mathcal T$. In fact, the graph polynomial specified previously and constructed from a rooted forest $\tilde{\mathcal F}_{e_{i_1} \cdots e_{i_k}}(\mathcal G)$, necessarily contain $\prod_{j=1}^k (e^{\beta H_{e_{i_1}}}-1)$. This means that if
 $\beta H_{e_{i_1}}=0$, necessarily these must be zero. Also, at high temperatures, this implies that  $\mathcal F_{e_1\cdots e_k}(\mathcal G)$ must be proportional to $\beta^k$. 
 As such, one plausibly choose a tree such that the orthogonal rooted forest are minimal, or that the tree is maximal, with graph edge variables $w_{ij}=\text{sup}|e^{\beta H_{ij}}-1|$.
 
 At this point, we can state the best choice of a reference tree.\ \\\ \\

 \textbf{Definition 2}. Let $\mathcal G$ be a sparse graph according to Definition 1 and $\beta\leq \beta^*$ be finite.  Let us define the weighted graph $\mathcal G_{w_{ij}}=\{w_{ij} A_{ij}\}$ according to the definition of $w_{ij}$ above. Given a tree $\mathcal T$ we define the stripped weighted graph $\bar {\mathcal G}_{w_{ij}}=\mathcal G_{w_{ij}}\setminus \mathcal T$ in which the elements of $\mathcal T$ are removed.
 Then, we say that a spanning tree is $\mathcal T$ optimal if
 \begin{eqnarray}
     \mathcal T^*=\text{min}_{\mathcal T} ||\bar { \mathcal G}_{w}||^2_2.
 \end{eqnarray} 
 where  $||\bar { \mathcal G}_{w}||^2_2$ is the Frobenius norm.
 It can be seen immediately that if $\mathcal G$ is a tree, such solution corresponds to the graph itself.\ \\\ \\
 
Let us now provide a way of evaluating these orthogonal sets via the tree $\mathcal T$. Let us now  define $\mathcal T_e=\mathcal T\setminus \tilde C_e$ as the portion of the reference spanning tree which does not contain $\tilde C_e$.

 Given an edge $e$, introduce the notation that $e\circ {\mathcal G}^\prime$ is the graph $G$ with the introduction of the edge $e$. Then we have\ \\

%\textbf{Lemma 4}
Let a graph $\mathcal G$ be a graph with a single cycle, and arbitrary spanning tree $\mathcal T$ and cotree $\bar {\mathcal T}= G\setminus \mathcal T$, we have
\begin{eqnarray}
    \mathcal F(e\circ  \mathcal T)&=&\mathcal F(\mathcal T)+ e\circ \big(\mathcal F(\mathcal T)-\mathcal F(\mathcal T_e)\circ \tilde C_e\big) \nonumber \\
    &=&\mathcal F(\mathcal T)+ e\circ \mathcal F(\mathcal T)-\mathcal F(\mathcal T_e)\circ  C_e
    \label{eq:comp}
\end{eqnarray}

%\textit{Proof}. 
We can prove the formula above as follows.
Any graph $\mathcal G$ can be constructed using the following procedure. Let $\mathcal T \cup \bar {\mathcal T} $ be a tree- cotree decomposition for $G$. Then let us call $e_i$, $i=1\cdots k$ be the edges of $\bar {\mathcal T}$. We have
\begin{eqnarray}
    \mathcal G=e_1\circ (e_2 \circ (\cdots (e_k \circ \mathcal T))).
\end{eqnarray}
Since $G$ contains only one cycle, we have that $|\bar {\mathcal T}|=1$.

Then, in order to prove eqn. (\ref{eq:comp}) it is sufficient to prove that, given a spanning tree $\mathcal T$, and $e\in \bar {\mathcal T}$
\begin{eqnarray}
    \mathcal F(e\circ \mathcal T)=\mathcal F(\mathcal T)+ e\circ \big(\mathcal F(\mathcal T)-\mathcal F(\mathcal T_e)\circ  \tilde C_e\big).
\end{eqnarray}
Now, given a spanning tree $\mathcal T$, $e\circ \mathcal T$ contains a cycle, and $\mathcal F(\mathcal T)$ contains $\mathcal T$ because it is a tree. Now it is easy to see that the only elements of the forest which contain the cycle are those elements of $\mathcal F(\mathcal T)$ which contain the reduced cycle $\tilde C_e$. These can be obtained as $\mathcal F(\mathcal T_e)\circ \tilde C_e$, which generates all the elements in $\mathcal F(\mathcal T)$ which contains $\mathcal C_e$.  $\square$

It follows that $Q(e_i)=\mathcal F(\mathcal T_e)\circ \tilde C_{e_i}$. This proves the formula above essentially for a graph close to a tree, e.g. when we have a graph which is a tree with  the addition of a single edge, or one cycle. When many cycles are present the formula becomes a little more involved, and one needs to introduce some new operations between cycles, e.g. a cycle calculus. It is very well known that the fundamental cycles form a basis for the cycle space. Given a set of fundamental cycles $\mathcal C$, we define the following two binary operations for $C_e,C_{e^\prime}\in \mathcal C$:

\begin{eqnarray}
    C_{e\circ e^\prime}=C_e\circ C_{e^\prime}&=&\{ e: e\in C_e \cup C_{e^\prime} \}\nonumber \\
   C_{e \cap e^\prime}=C_e \cap C_{e^\prime}&=& \{ e: e\in C_e \cap C_{e^\prime} \} \nonumber \\
    \tilde C_{e\square e^\prime}=C_e\square C_{e^\prime}&=&\{ e: e\in (C_e \cup C_{e^\prime})\setminus (C_e \cap C_{e^\prime}) \}\nonumber \\
\end{eqnarray}
and similar notions for $\tilde C_e$ and $\tilde C_{e^\prime}$. Note that if $(C_e \cap C_{e^\prime})=\{\emptyset\}$, $C_{e\circ e^\prime}=C_{e\square e^\prime}$, and that we can write  $C_{e\circ e^\prime}=e\circ e^\prime \circ \mathcal T$.

On the other hand $ C_{e\square e^\prime}$ is in fact a new cycle contained in $\mathcal G$ and obtained from the cycle basis. Thus, similarly to the case of a single cycle, if $e_1,e_2\in \bar {\mathcal T} $,  $e_1\circ e_2 \circ \mathcal F(\mathcal T)$
adds to the forest all graphs which contain both $e_1$ and $e_2$. Naturally, however, these are not all trees, and thus some of these must be subtracted.
As in the case of a single cycle, such cycle is local to $C_{e_1}$ and $C_{e_2}$ (and this is because $C_e$ is a basis), and thus the graphs that need to be subtracted are all those that contain some cycles.
These are all graphs that contain as a subgraph $C_{e_1}$ and/or $C_{e_2}$, or $ C_{e_1 \square e_2}$ if $C_{e_1 \square e_2} \neq C_{e_1 \circ e_2}$.
These subgraphs are 
\begin{eqnarray}
 & & e_1\circ e_2 \circ \tilde C_{e_1\circ e_2},\ \ \ e_1 \circ e_2 \circ \tilde C_{e_1},\nonumber \\
& & e_1 \circ e_2 \circ \tilde C_{e_2},\ \ \  e_1 \circ e_2 \circ \tilde C_{e_1\square e_2}.
\end{eqnarray}
Let us call 
\begin{eqnarray}
Q(e_1,e_2)=e_1\circ e_2 \circ (\tilde C_{e_1\circ e_2}&+&  \tilde C_{e_1}+   \tilde C_{e_2} \nonumber \\
&+&(1-\delta_{e_1 \circ e_2,e_1\square e_2})  \tilde C_{e_1\square e_2}).\nonumber 
\end{eqnarray}

In the equation above, $\delta_{\delta_{e_1 \circ e_2,e_1\square e_2}}$ is a shorthand notation; $\delta_{\delta_{e_1 \circ e_2,e_1\square e_2}}=1$ if $C_{e\circ e^\prime}=C_{e\square e^\prime}$ and zero otherwise.
We defined a reduced tree as $\mathcal T_{ee^\prime}=\mathcal T\setminus  C_{e \square e^\prime}$, which is the tree which does not contain elements of $\tilde C_{e_1 \cup e_2}=\tilde C_{e_1} \cup \tilde C_{e_2}$. Then, similarly to the case of a single cycle, we have for two adjacent cycles, we have
\begin{eqnarray}
    \mathcal F(e_1 \circ e_2 \circ \mathcal T)&=&(1+e_1 + e_2 +e_1 \circ e_2)\circ \mathcal F(\mathcal T) \nonumber \\
    &-&\mathcal F(\mathcal T_{e_1}) \circ C_{e_1}-\mathcal F(\mathcal T_{e_2}) \circ C_{e_2} \nonumber \\
    &-&\mathcal F(\mathcal T_{e_1 e_2}) \circ Q(e_1,e_2).
\end{eqnarray}
At this point however we should start to see a pattern. If we introduce a third cycle, we will have $e_1\circ e_2\circ e_3 \circ \mathcal F(\mathcal T)$ but we should remove from these the one, two and three cycles contribution. These are 
\begin{widetext}
\begin{eqnarray}
    \mathcal F(e_1 \circ e_2 \circ e_3 \circ \mathcal T)&=&(1+e_1 + e_2 +e_3+e_1\circ e_2+e_2\circ e_3+e_1\circ e_3+e_1 \circ e_2\circ e_3)\circ \mathcal F(\mathcal T) \nonumber \\
    &-&\mathcal F(\mathcal T_{e_1}) \circ C_{e_1}-\mathcal F(\mathcal T_{e_2}) \circ C_{e_2}-\mathcal F(\mathcal T_{e_3}) \circ C_{e_3} \nonumber \\
    &-&\mathcal F(\mathcal T_{e_1 e_2}) \circ Q(e_1,e_2)-\mathcal F(\mathcal T_{e_2 e_3}) \circ Q(e_2,e_3)-\mathcal F(\mathcal T_{e_1 e_3}) \circ Q(e_1,e_3) \nonumber \\
    &-& \mathcal F(\mathcal T_{e_1 e_2 e_3}) \circ Q(e_1,e_2,e_3)\nonumber \\
    Q(e_1,e_2,e_3)&=&e_1 \circ e_2 \circ e_3 \circ (\tilde C_{e_1\circ e_2 \circ e_3}+  \tilde C_{e_1\circ e_2} + \tilde C_{e_1\circ e_3}+  \tilde C_{e_2\circ e_3} \nonumber \\
    &+& (1-\delta_{\circ,\square})\tilde C_{e_1\square e_2} + (1-\delta_{\circ,\square})\tilde C_{e_1\square e_3}+  (1-\delta_{\circ,\square})\tilde C_{e_2\square e_3}+(1-\delta_{\circ,\square})C_{e_1\square e_2\square e_3})
\end{eqnarray}
\end{widetext}
where we used a shorthand notation for the delta, and is the expression for the 3rd order correction.
We now note something interesting, which is that what the calculation above suggest is the introduction of a ``loop" algebra, as advocated in \cite{chertkov,chertkov2}, and we thus re-obtain (with a different technique), the corrections obtained via belief propagation for the case of graphical models.

Let us now note that the cycle operator acting on $\mathcal F(\mathcal T)$ can be written, recursively, as 
\begin{eqnarray}
    \prod_{e_i\in \bar{\mathcal T}}(1+e_i) \mathcal F(\mathcal T)&=&e^{-\beta H}|_{\bar {\mathcal T}} \mathcal F(\mathcal T) \nonumber \\
    &=&e^{-\beta H}|_{\bar {\mathcal T}} e^{-\beta H}|_{\mathcal T}= e^{-\beta H}.
\end{eqnarray}
To our surprise we are thus reconstructing, cycle by cycle, the full probability distribution. For computability purposes reconstructing the full forest expansion might is as hard as computing the full model. What have we gained this far?  We could have obtained such an approximation directly from $e^{-\beta H}$, decomposing $\mathcal G=(\mathcal T,\bar {\mathcal T})$ as
\begin{eqnarray}
    e^{-\beta H}&=&e^{\beta \sum_{(ij)\in \mathcal G} A_{ij} H_{ij}}=e^{\beta \sum_{(ij)\in \mathcal T} A_{ij} H_{ij}} e^{\beta \sum_{(ij)\in {\bar {\mathcal T}} A_{ij} H_{ij}}} \nonumber \\
    &=&e^{\beta \sum_{(ij)\in \mathcal T} A_{ij} H_{ij}}\prod_{(ij)\in \bar {\mathcal T}} (1+(e^{\beta H_{ij}}-1)) \nonumber \\
    &=&(1+\cdots )e^{-\beta H}|_{\mathcal T}
\end{eqnarray}

However, we would have missed the corrections to the tree approximation required by the forest expansion:  we are iteratively  removing the cycles from the decomposition in order to regain the right sum.

\begin{figure}
    \centering
    \includegraphics{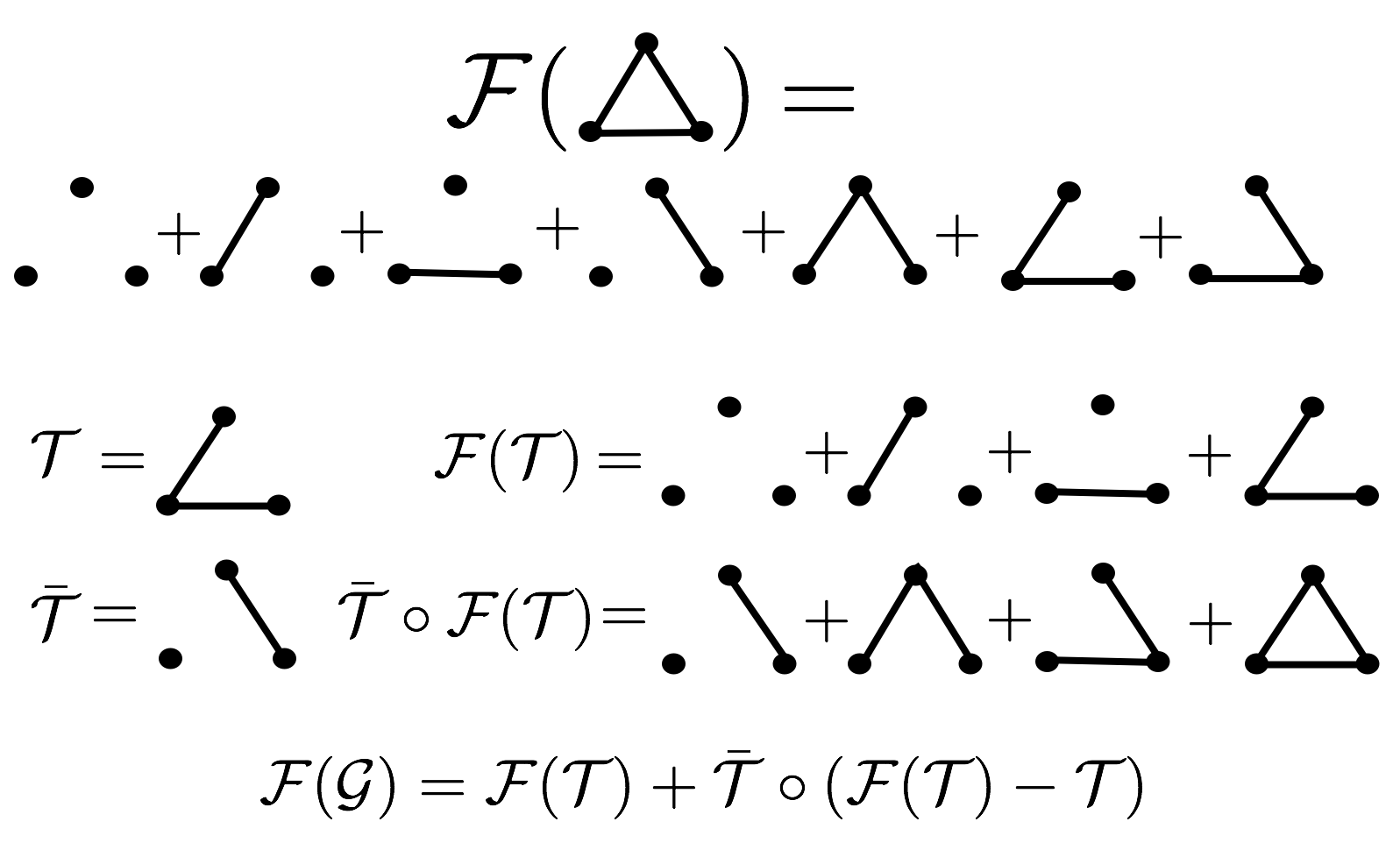}
    \caption{The tree and co-tree decomposition of the forest expansion of $K_3$, in which the decomposition of eqn. (\ref{eq:comp}) is particularly simple. }
    \label{fig:treecotree}
\end{figure}

A particular simple case is the one of $\mathcal G=K_3$, this is shown in Fig. \ref{fig:treecotree}. There we have that $\tilde C_e=\mathcal T$ and $e=\bar {\mathcal T}$, and $\mathcal T_e=\{\emptyset\}$ and $\mathcal F(\mathcal T_{e})=1$. Thus $ \mathcal F(K_3)=\mathcal F(\mathcal T)+\bar {\mathcal T} \circ (\mathcal F(\mathcal T)-\mathcal T)$ .

\subsection{First order approximation}

Note that eqn. (\ref{eq:comp}) explicitly shows how the cycle decomposition enters in the forest formula for an arbitrary graph. Clearly, we have an immense freedom in the choice of $\mathcal T$, but this is not uncommon, as it happens also for resistive circuits.

Let us see how the first order approximation to the partition function looks like. Following the discussion of the previous section, given the graph $\mathcal G$ and the edge weights $w_{ij} A_{ij}=\beta \text{sup} |H_{ij}|$,  we pick a maximally spanning tree $\mathcal T$. Then, we know that
\begin{eqnarray}
    e^{-\beta H}&\approx& \mathcal F(\mathcal T)+\sum_{i=1}^K e_i \circ (\mathcal F(\mathcal T)-\mathcal F(\mathcal T_{e_i})\circ \tilde C_{e_i})+\cdots.
\end{eqnarray}

Let us now see what formula (\ref{eq:comp}) implies for $Z_{\bar {\mathcal F}}$.  The operator $\circ$ is simply the multiplication in terms of the polynomials. 
The advantage of writing the expression in terms of $\mathcal F(\mathcal T)$ and $\mathcal F(\mathcal T_e)$ is that these immediately map to Boltzmann distributions, via an intermediate summation over the forest, on the tree.
We thus have the mapping
\begin{eqnarray}
    \mathcal F(\mathcal T)&\rightarrow& \prod_{(ij)\in \mathcal T}e^{\beta H_{ij} } \nonumber \\
   e_{kt}\circ \mathcal F(\mathcal T)&\rightarrow& (e^{\beta H_{kt}}-1) \prod_{(ij)\in \mathcal T}e^{\beta H_{ij} } \nonumber \\
    \mathcal F(\mathcal T_e)&\rightarrow& \prod_{(ij)\in \mathcal T_e}e^{\beta H_{ij} } \nonumber \\
    C_e&\rightarrow&\prod_{(ij)\in C_e} (e^{\beta H_{ij}}-1)
\end{eqnarray}
because of the expansion we have explained at the beginning of this section.

As a simple example of the expansion above, let us consider the limits
\begin{eqnarray}
    \lim_{T\rightarrow \infty}\lim_{L\rightarrow \infty} Z_{\bar {\mathcal F}} \approx Z_{\mathcal T}.
\end{eqnarray} 
Note that swapping the limits is a highly non-trivial statement.
In the high temperature approximation can write 
\begin{eqnarray}
    e^{\beta H_{ij}}-1\approx \beta H_{ij}+O(\beta^2).
\end{eqnarray}
It follows that, at the first order of the expansion, we have a high temperature expansion of the form
\begin{eqnarray}
    e^{-\beta H}|_{L\rightarrow \infty}&\approx& (1+\beta \sum_{e_{kt} \in \bar {\mathcal T}} H_{kt}) e^{-\beta H}|_{\mathcal T} \nonumber \\
    &-&  \sum_{e\in \bar {\mathcal T}} \beta^{|C_e|}(\prod_{e_{jt}\in C_e} H_{jt}) e^{-\beta H}|_{\mathcal T_e}
    \label{eq:expansionbeta}
\end{eqnarray}
which is one of the main results of this paper. 
The expansion above is in fact exact for the case of $1$-dimensional models on a ring of size $L\gg 1$, and the summation is reduced to only one element $e$ (the weakest link in the chain if it exist). 
Something to note is that both $e^{-\beta H}|_{\mathcal T}$ and $e^{-\beta H}|_{\mathcal T_e}$ are being evaluated on trees.

\section{Application}
\subsection{Examples}
First, we note that for a one-dimensional model on a circle of length $L$, the forest expansion explicitly says that the ``open" partition function is almost identical to the ``closed" one. This is because $\Delta F=0$ identically in this case. This example is thus a little too simple.

\begin{figure}
    \centering
    \includegraphics{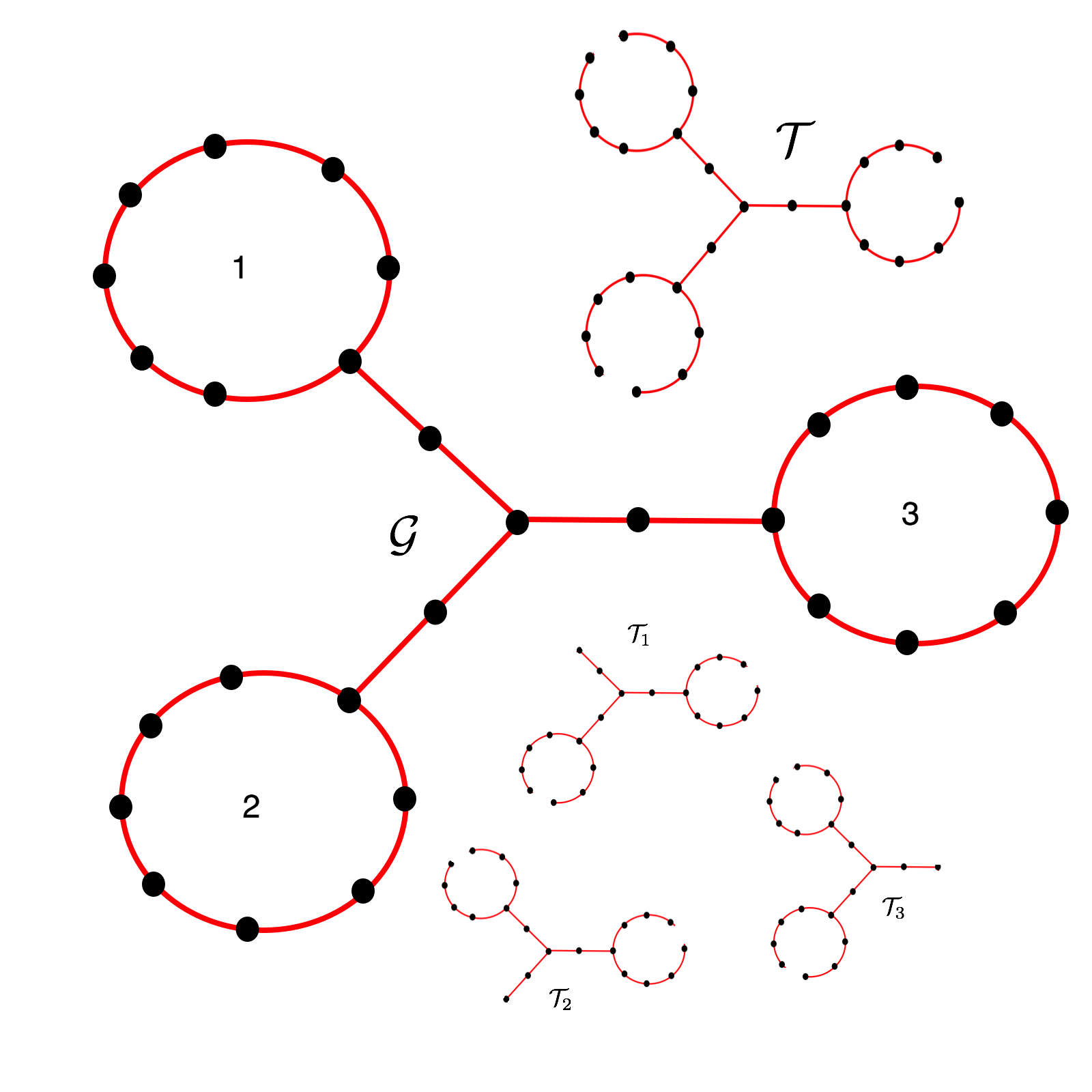}
    \caption{An example of the tree approximation to a ``loopy" graph, and the cycle reduced trees that occur at the first order approximation.}
    \label{fig:example}
\end{figure}

Let us thus provide some examples of reduced trees in the case of an interaction graph like in Fig. \ref{fig:example}. According to eqn. (\ref{eq:expansionbeta}), given $\mathcal T$, we need to compute the cycles, and the reduced trees associated to these. In  \ref{fig:example} we have three cycles. Thus, we can write
\begin{eqnarray}
    \mathcal F(\mathcal G)=\mathcal F(\mathcal T)+\Delta F
\end{eqnarray}
In that case, we have $\Delta F$ is the sum over three corrections. 
Let us call $e_1,e_2,e_3$ the three edges that are removed from $\mathcal G$ to obtain $\mathcal T$. Then, 
\begin{eqnarray}
    \Delta F_1=\beta H_{e_1} e^{-\beta H}|_{\mathcal T}-e^{-\beta H}|_{\mathcal T_1}\prod_{e\in \mathcal C_1} \beta H_{e}.
\end{eqnarray}
and similar expressions for $\Delta F_2$ and $\Delta F_3$.

Another common example of two-body interaction is the Ising model, defined by:
\begin{equation}
    Z(\beta, A )=\sum_{\{\sigma=\pm 1\}} e^{\beta  \sum_{ij} \sigma_i \sigma_j J_{ij}},
\end{equation}
where we assume $J_{ij}= J A_{ij}$. It is known that evaluating such partition function is an NP-hard problem if $A_{ij}$ is not a tree. 

First, let us rephrase eqn. (\ref{eq:q1}) for the case of the Ising model. In fact, we can write for binary variables $\sigma_i=\pm 1$
\begin{eqnarray}
    e^{A \sigma_i \sigma_j}&=&\cosh( A \sigma_i \sigma_j)+\sinh(A\sigma_i \sigma_j) \nonumber \\
    &=&\cosh( A )+\sinh(A)\sigma_i \sigma_j
\end{eqnarray}
It follows that $\cosh( H_{ij} \sigma_i \sigma_j \text{min}(u_1,\cdots,u_{L-1}))$ and $\sinh( H_{ij} \sigma_i \sigma_j \text{min}(u_1,\cdots,u_{L-1}))$ contribute to the even and odd series expansion for $e^{H_{ij} \sigma_i \sigma_j \text{min}(u_1,\cdots,u_{L-1})}$.
Something very interesting which occurs for binary variables $\sigma_i =\{\pm 1\}$ is however the following. Let $C_{e}$ be a cycle in $\mathcal G$. Given the fundamental cycles $C_e$, we define the closed cycle observables given by
\begin{eqnarray}
    \mathcal L(e)=\prod_{e^\prime \in \mathcal C_e} H_{e^\prime}=\prod_{v_j \in C_e}^{|C_e|} \beta H_{v_j v_{j+1}}.
\end{eqnarray}
For the Ising models, for a closed cycle of length 3 this is given by quantities of the form
\begin{eqnarray}
    \mathcal L_3(e)=\beta^3 H_{v_1 v_2}H_{v_2 v_3}H_{v_3 v_1}.
\end{eqnarray}
Now note that quantities of this form are independent from the variables $\sigma$ for the Ising model, as $\sigma_i^2=1$. As a consequence, for the Ising model we obtain a first order expansion correction due to eqn. (\ref{eq:expansionbeta}), of the form
\begin{eqnarray}
    e^{-\beta H}&\approx& (1+\beta \sum_{k=1}^{|\bar {\mathcal T}|} H_{v_{(e_k)_1} v_{(e_k)_2}} \sigma_{v_{(e_k)_1}} \sigma_{ v_{(e_k)_2}})  e^{-\beta H}|_{\mathcal T} \nonumber \\
    & &\ -\sum_{k=1}^{|\bar {\mathcal T}|} c_k e^{-\beta H}|_{\mathcal T_k} 
\end{eqnarray}
where $\mathcal T$ is a spanning tree, $\mathcal T_k$ are reduced trees in which each element of $\mathcal T$ which is contained in $\tilde C_{e_k}$ has been removed, and $c_k$  are constants.

%The next example is the simplest possible model, e.g. a partition function of the form
%\begin{eqnarray}
%    \mathcal Z=\int_{-\infty}^\infty d x_i e^{-\frac{1}{2} \sum_{ij} x_i A_{ij} x_j}=\frac{\sqrt{(2\pi)^N}}{\sqrt{\text{det}(A)}}
%\end{eqnarray}
%\textcolor{red}{consider the case of $A$= spanning tree+edge}
\subsection{Numerical study using belief propagation}
In order to investigate the temperature at which the convergence to the tree expansion applies, we studied the Ising model $H=-\frac{J}{2} \sum_{ij} A_{ij} \sigma_i \sigma_j$, with homogeneous coupling  on trees without external field, where $A_{ij}$ is the undirected adjacency matrix of the graph. Specifically, we generated random spanning trees $\mathcal T$ on $N$ nodes (we first used a uniform random matrix $\tilde A$, and then applied a minimum spanning tree Kruskal algorithm), and then added an edge between two nodes in order to introduce a loop, generating a graph $\mathcal T^\prime$. Typically, such loop will be large and thus we should expect loopy belief propagation to converge, and a (maximum spanning) tree approximation to work. This is well known in the literature. 

\begin{figure}
    \centering
    \includegraphics[scale=0.18]{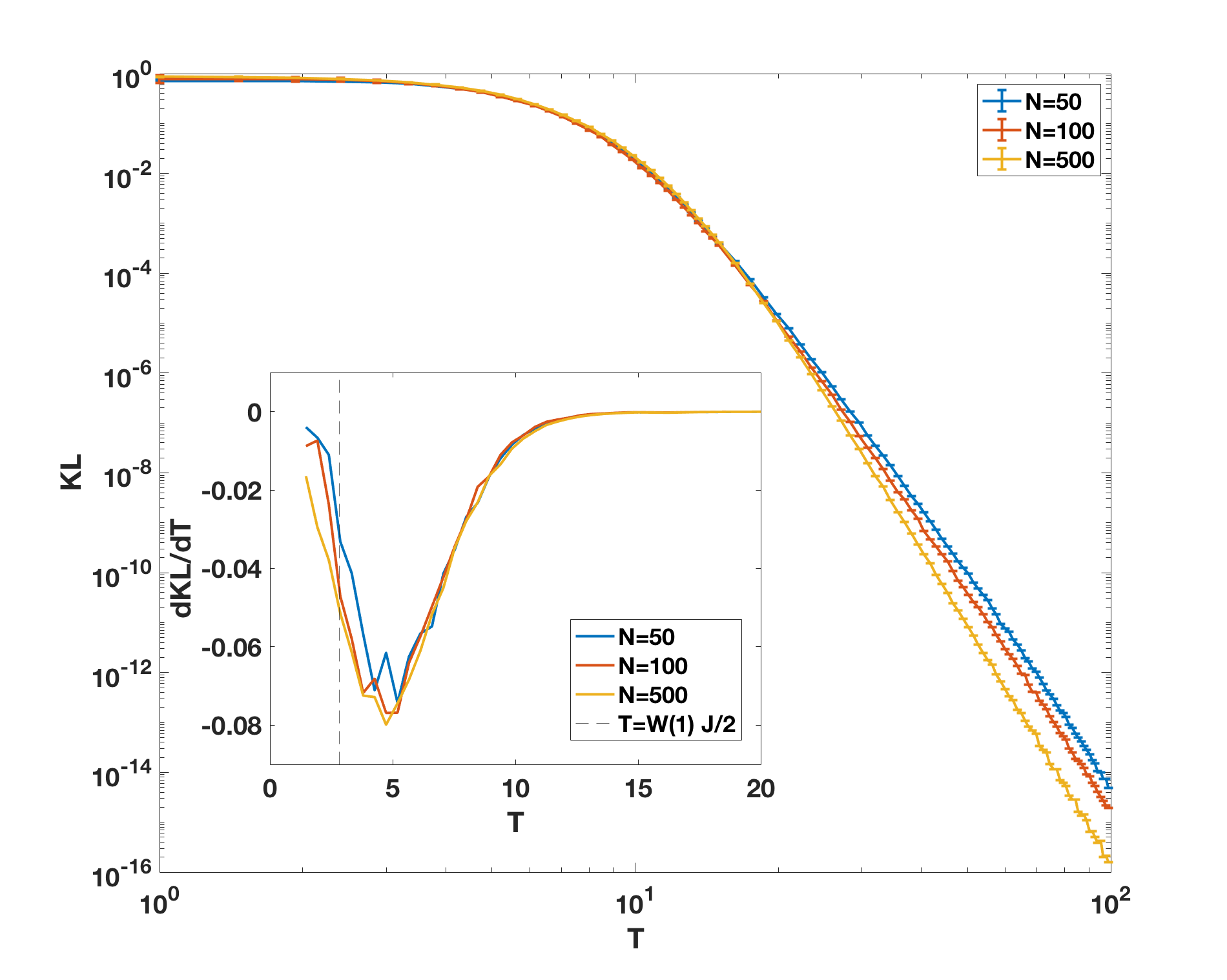}
    \caption{Average Kullback-Leibler divergence averaged over 1000 different trees between the beliefs calculated using the exact belief propagation for $J=10$, and for tree sizes $N=50,100,500$. The results are consistent among different sizes, and we see a decline of the KL divergence for $T\approx W(1) J$.}
    \label{fig:kl}
\end{figure}

We then used belief propagation in order to obtain the beliefs $p(s_i)$ for each spin $s_i$, both in the case of the loopy tree $\mathcal T^\prime$ and $\mathcal T$. In order to measure how close the beliefs were, we used the  Kullback-Leibler divergence between two distributions, $KL(P||Q)=\sum_{\{s_i=\pm 1\}} p(s_i) \log \frac{p(s_i)}{q(s_i)}$. We then calculated the average over every spin, 
$KL(P_{\mathcal T}||P_{\mathcal T^\prime})=\frac{1}{N} \sum_{i=1}^N KL(p_{\mathcal T}(s_i)||p_{\mathcal T^\prime}(s_i))$. Clearly, we have $KL(P_{\mathcal T}||P_{\mathcal T^\prime})=0$ if the two beliefs are identical. We have further averaged the KL divergence over $N_{MC}=1000$ samples, for tree sizes of $N=50,100,500$.  The results are shown in Fig. \ref{fig:kl}.  According to the prediction of this article, we should observe that for temperatures above $T\approx 0.567 J/2$ the two provide the same result, e.g. $P_{\mathcal T}\rightarrow P_{\mathcal T^\prime}$. We studied both $KL$ and $\frac{dKL}{dT}$ (inset), as a function of the temperature for $T\in [0.1, 100]$. For $J=10$, we should expect that at $T\approx 5.67/2$ a transition between the two regimes occurs. This is intuitively what we see in Fig. \ref{fig:kl}. One should expect such transition to occur at $T> J/2$, but indeed it occurs at an earlier temperature.
\section{Conclusions}
In this paper we have put forward a methodology to write the Boltzmann probability distribution of a 2-body Hamiltonian in terms of trees, following a formula first derived by Brydges and Kennedy, and then refined by Abdessalam and Rivasseau (the BKAR expansion). Such formula is a generalization of the fundamental theorem of calculus to functions of $N(N-1)/2$ variables and the expansion for the partition function can be thought of as a generalization of polymeric expansions.  While re-deriving the loop expansion first obtained by Chertkov and Chernyak, we have shown that the ``cycle perturbative" phase is hidden within a ``forest phase", with two different expansion temperatures. First, we have shown that 
for sparse graphs with large cycles (and bounded interactions), the probability distribution can be written in terms of a polynomial forest expansion,e.g. sum over connected or disconnected subtrees of the interaction graph, for a temperature such that $\beta\leq \beta^*\approx 0.54 \text{sup}_{ij}|H_{ij}|$. 
We have in particular focused on two regimes: graphs with large cycles which are cycle dense or cycle sparse. Our definition of dense is based on the scaling of the degree of the dual graph (e.g. cycles become vertices). If the degree of the dual graph scales with the size of the interaction network, then the graph is dense according to this definition. The difference between the two cases is only important when $\beta=\beta^*$, in which case for dense graphs the correction to the forest sum is of order one, whereas for sparse graphs if it is of order $1/L$, with $L$ the minimum cycle length.

We have then shown that we can recast the forest expansion in terms of a tree-cotree splitting of the original interaction Hamiltonian. We have also provided a favorable choice of such tree, e.g. a maximally spanning tree on a weighted interaction graph.

The cycle perturbative phase result should not come as a surprise. For quite some time we have known that for sparse graphs belief propagation (which converges only on trees), provides reasonable approximation to the probability distributions. In this sense, in this paper we provide an alternative venue and background to this statement.
However, the fact that this phase is hidden in a ``forest phase" is a novel result of this article.

The techniques we used in the paper have however some important practical implications, of which can be read directly from our constructive proof. For instance, many models on trees can be solved exactly (at least numerically) via techniques such as cavity methods or message passing. In this sense, this paper establishes that there is a particular temperature $\beta^*$ above which models on (large) loopy graphs can be solved by approximating the interaction using an optimal choice of spanning tree. Temperature corrections can also be calculated at higher orders both in temperature and cycles. However, we have shown that these can be also cast in terms of reduced trees, and thus are also \textit{easy} to calculate via belief propagation.
We tested the fact that $T^*=W(1) J$ is a special temperature for a tree approximation with a numerical experiment for the Ising model and with the use of belief propagation. While one would expect a tree approximation to work for $T>J$, we observed numerically that the KL divergence (which we use as a proxy for the tree approximation), starts to decline at $T^*\approx J/2$, which is consistent with $T^*_{th}=  W(1)J$.

As a final comment, one might wonder why this approach might be problematic for a quantum system. The reason is that when one takes the derivatives of the probability distribution (which is an operator) and defined in terms of $\hat H=\sum_{ij} \hat H_{ij}$, not necessarily one has that $[\hat H_{ij},\hat H_{i^\prime j^\prime}]$ when $i=i^\prime$ or $j=j^\prime$. For this reason, while possible, extensions of the result of the present work to quantum systems are out of the scope of this paper and for future investigations.
\ \\\ \\

\textit{Acknowledgments.} We would like to thank A. Lokhov and M. Vuffray for various discussions. We acknowledge the support of NNSA for the U.S. DoE at LANL under Contract No. DE-AC52-06NA25396, and via DOE-LDRD grants PRD20170660 and PRD20190195.  
No trees were cut or harmed during the preparation of this paper.

\end{document}